\documentclass[11pt]{article}

\usepackage{latexsym,amsmath,amsfonts,amssymb,amstext,theorem,euscript,array,enumerate,mathrsfs}
\usepackage{color}
\usepackage{comment}

\newtheorem{Theorem}{Theorem}[section]
\newtheorem{Definition}{Definition}[section]
\newtheorem{Proposition}{Proposition}[section]

\newtheorem{Lemma}{Lemma}[section]
\newtheorem{Corollary}{Corollary}[section]
\newtheorem{Remark}{Remark}[section]

\numberwithin{equation}{section}

\def\esssup_#1{\underset{#1}{\mathrm{ess\,sup\, }}}
\def\essinf_#1{\underset{#1}{\mathrm{ess\,inf\, }}}

\def \trans{^{\scriptscriptstyle{\intercal}}}

\def\sqr#1#2{{\vcenter{\vbox{\hrule height .#2pt \hbox{\vrule
 width .#2pt height#1pt \kern#1pt \vrule
width .#2pt} \hrule height .#2pt}}}}

\def\ds{\begin{displaystyle}}
\def\eds{\end{displaystyle}}

\def\<{\langle }
\def\>{\rangle }

\def \Inf{\displaystyle\inf}
\def \Sup{\displaystyle\sup}

\def \N{\mathbb{N}}
\def \R{\mathbb{R}}
\def \M{\mathbb{M}}

\def \E{\mathbb{E}}
\def \F{\mathbb{F}}

\def \P{\mathbb{P}}
\def \Q{\mathbb{Q}}

\def\a{{\bf a}}
\def\b{{\bf b}}
\def\p{{\bf p}}
\def\M{{\bf M}}

\def \Ac{{\cal A}}
\def \Bc{{\cal B}}
\def \Cc{{\cal C}}
\def \Dc{{\cal D}}
\def \Ec{{\cal E}}
\def \Fc{{\cal F}}
\def \Gc{{\cal G}}
\def \Hc{{\cal H}}

\def \Kc{{\cal K}}

\def \Nc{{\cal N}}

\def \Wc{{\cal W}}

\def \eps{\varepsilon}

\def \ep{\hbox{ }\hfill$\Box$}
\def \epR{\hbox{ }\hfill$\lozenge$}

\def\Dt#1{\frac{\partial #1}{\partial t}}

\def\reff#1{{\rm(\ref{#1})}}

\def\beqs{\begin{eqnarray*}}
\def\enqs{\end{eqnarray*}}
\def\beq{\begin{eqnarray}}
\def\enq{\end{eqnarray}}

\allowdisplaybreaks

\begin{document}

\title{Zero-sum stochastic differential games of \\ generalized McKean-Vlasov type\thanks{This research was initiated when A. Cosso visited LPSM under the support of the ANR project CAESARS (ANR-15-CE05-0024).
}}

\author{
Andrea COSSO\footnote{Department of Mathematics, University of Bologna, Piazza di Porta S. Donato 5, 40126 Bologna, Italy, \sf andrea.cosso at unibo.it}
\qquad\quad
Huy\^en PHAM\footnote{LPSM, University Paris Diderot and CREST-ENSAE, \sf pham at math.univ-paris-diderot.fr}
}

\maketitle

\begin{abstract}
We study zero-sum stochastic differential games where the state dynamics of the two players is governed by a generalized McKean-Vlasov (or mean-field) stochastic differential equation in which the distribution of both state and controls of each player appears in the drift and diffusion coefficients, as well as in the running and terminal payoff functions.  We prove the dynamic programming principle (DPP) in this general setting, which also includes the control case with only one player, where it is the first time that DPP is proved for open-loop controls. We also show that the upper and lower value functions are viscosity solutions to a corresponding upper and lower Master Bellman-Isaacs equation. Our results extend  the seminal work of Fleming and Souganidis \cite{FlemingSouganidis} to the McKean-Vlasov setting.
\end{abstract}

\vspace{7mm}


\vspace{7mm}


\date{}

\vspace{7mm}

\section{Introduction}

McKean-Vlasov (McKV) control problem (also called mean-field type control problem) has been knowing a surge of interest with the emergence of mean-field game (MFG) theory, see \cite{lio12}, \cite{benetal13}, \cite{phawei17}, \cite{cardelbook}. 
 Such a problem was originally motivated by large population stochastic control under mean-field interaction in the limiting case where the number of agents tends to infinity; now various applications can be found in economics, finance, and also in social sciences for modeling motion of socially interacting individuals and herd behavior. A crucial assumption in large population models, as well as in the theory of MFG and McKV control problems, is the homogeneity of agents with identical outcomes.

In this paper, we are concerned with generalized McKean-Vlasov 
stochastic differential equations controlled by two players with opposite objectives: this problem is then called zero-sum stochastic differential game of generalized McKean-Vlasov type. A typical motivation arises from the consideration of two infinite groups of homogenous interacting agents pursuing conflicting interests that can also interact between each other, like e.g. in pursuit/evasion games. 
 
 The seminal paper \cite{FlemingSouganidis} formulated in a rigorous manner zero-sum stochastic differential games, with state dynamics governed by standard stochastic differential equations. The formulation in \cite{FlemingSouganidis} can be described as {\it nonanticipative strategies against open-loop controls}. This distinction of controls type between the two players  is crucial to show the dynamic programming principle (DPP) for the upper and lower value functions (which coincide, i.e. the game has a value, under the so-called Isaacs condition), as it is known that the formulation {\it open-loop controls vs open-loop controls} does not give rise in general to a dynamic game and a fortiori to a DPP, see for instance Buckdahn's counterexample in Appendix E of \cite{phazha14}.  
 
 Zero-sum McKV stochastic differential games were recently considered in \cite{limin16} and \cite{djeham16} in a weak formulation where only the drift (but not the diffusion coefficient) depends on controls and state distribution. Notice that as the authors work on a canonical probabilistic setting, their game can be seen as a game in the form {\it feedback controls vs feedback controls}.  
 We mention also the recent paper \cite{aver18}, which considers deterministic mean-field type differential games with feedback controls. 
 
In the present work, we study zero-sum stochastic differential games of generalized McKean-Vlasov type where all the coefficients of both state dynamics and payoff functional depend upon the distributions of state and controls (actually, they can also depend on the \emph{joint} distribution of state and controls, however under the standard continuity and Lipschitz assumptions, it turns out that the coefficients only depend on the marginal distributions, see Remark \ref{R:JointLaw}). As in \cite{FlemingSouganidis}, we use a strong formulation with nonanticipative strategies against open-loop controls. We define the lower and upper value functions of this game, and our first contribution (Proposition \ref{P:Lift}) is to show that they can be considered as functions on the Wasserstein space of probability measures. Notice that this is a nontrivial issue as we do not restrict to  feedback controls.
 
Our second main result (Theorem \ref{T:DPP}) is the proof of the dynamic programming principle  for the lower and upper value functions. The key observation is to reformulate the problem as a deterministic differential game in the infinite dimensional space $L^q$ of $q$-th integrable random variables. Notice that the proof of the DPP is relevant also for the control case (which corres\-ponds to the special case where the space of control actions of the second player is a singleton), as a matter of fact in the present paper we consider open-loop controls, while in the literature the DPP has been proved only for feedback controls,  
see \cite{phawei17}. Let us mention however the paper \cite{bayetal} which states a \emph{randomized} DPP for the control case  with open-loop controls (but without dependence on the control distribution).   
We also show how to recover the standard DPP in the case without mean-field dependence.

The third contribution of this paper (Theorem \ref{prop:visco}) is the partial differential equation characterization of the value functions. By relying on the notion of differentiability in the Wasserstein space due to P.L. Lions, we prove the viscosity property of the lower and upper value functions to the corresponding dynamic programming  lower and upper Bellman-Isaacs equations. Uniqueness is stated when working on the lifted  Hilbert space $L^2$ 
and consequently, existence of a game value is obtained under a generalized Isaacs  condition. 

The outline of the paper is as follows. Section \ref{secformul} formulates the zero-sum stochastic diffe\-rential game of generalized McKean-Vlasov type. 
In  Section \ref{secprop} we show that the upper and lower value functions can be defined as functions on the Wasserstein space of probability measures.  
Section \ref{secDPP} is devoted to the rigorous statement and proof of the dynamic programming principle for both value functions. Finally, in Section \ref{secvisco} we prove the viscosity property of the value functions.

\section{Formulation of the stochastic differential game} \label{secformul}

Let $(\Omega,\Fc,\P)$ be a complete probability space on which a $d$-dimensional Brownian motion $W=(W_t)_{t\geq0}$ is defined. Let $\F^o=(\Fc_t^o)_{t\geq 0}$ be the filtration generated by $W$, and let $\F=(\Fc_t)_{t\geq0}$ be the augmentation of $\F^o$ with the family $\Nc$ of $\P$-null sets of $\Fc$, so that $\Fc_t=\Fc_t^o\vee\Nc$, for every $t\geq0$. Notice that $\F$ satisfies the usual conditions of $\P$-completeness and right-continuity. We also define, for every $t\geq0$, the filtration $\F^t=(\Fc_s^t)_{s\geq t}$ which is the $\P$-completion of the filtration generated by the Brownian increments $(W_s-W_t)_{s\geq t}$ (notice that, when $t=0$, $\F^0$ coincides with $\F$). Finally, we suppose that there exists a sub-$\sigma$-algebra $\Gc$ of $\Fc$ which is independent of $\Fc_\infty$, which will be assumed ``rich enough'', as explained below.

We fix a positive integer $n$ and a real number $q\in[1,\infty)$. We denote by $\mathscr P_{\text{\tiny$q$}}(\R^n)$ the family of all probability measures $\mu$ on $(\R^n,\Bc(\R^n))$ with finite $q^{\text{th}}$-order moment, i.e. $\|\mu\|_q$ $:=$ $\big(\int |x|^q \mu(dx)\big)^{1/q}$ $<$ $\infty$.  
More generally, for any $r\in[1,\infty)$ and $m\in\N\backslash\{0\}$ we denote by $\mathscr P_{\text{\tiny$r$}}(\R^m)$ the family of all probability measures on $(\R^m,\Bc(\R^m))$ with finite $r^{\text{th}}$-order moment. We endow 
$\mathscr P_{\text{\tiny$r$}}(\R^n)$ with the topology induced by the Wasserstein metric of order $r$:
\[
\Wc_{\text{\tiny$r$}}(\mu,\mu') = \inf\bigg\{\bigg(\int_{\R^n\times\R^n} |x - x'|^r\,\boldsymbol\mu(dx,dx')\bigg)^{1/r}\colon\boldsymbol\mu\in\mathscr P_{\text{\tiny$r$}}(\R^n\times\R^n)\text{ with marginals $\mu$ and $\mu'$}\bigg\},
\]
for all $\mu,\mu'\in\mathscr P_{\text{\tiny$r$}}(\R^n)$. We assume that the sub-$\sigma$-algebra $\Gc$ is ``rich enough'' in the following sense: $\Gc$ satisfies
\begin{equation}\label{G}
\mathscr P_{\text{\tiny$1$}}(\R) \ = \ \big\{\P_{\text{\tiny$\xi$}}\colon\xi\in L^1(\Omega,\Gc,\P;\R)\big\}, \qquad\qquad \mathscr P_{\text{\tiny$q$}}(\R^n) \ = \ \big\{\P_{\text{\tiny$\xi$}}\colon\xi\in L^q(\Omega,\Gc,\P;\R^n)\big\},
\end{equation}
where $\P_{\text{\tiny$\xi$}}$ denotes the distribution of $\xi$. Possibly making $\Gc$ smaller, we can suppose that there exists a random variable $\Gamma^\Gc\colon(\Omega,\Gc)\rightarrow(G,\mathscr G)$, taking values in some Polish space $G$ with Borel $\sigma$-algebra $\mathscr G$, such that $\Gamma^\Gc$ has an atomless distribution and $\Gc=\sigma(\Gamma^\Gc)$ (see Remark \ref{R:G}).

\begin{Remark}\label{R:G}
{\rm
It is well-known (see e.g. Theorem 3.19 in \cite{Kallenberg}) that the probability space $([0,1],\Bc([0,1]),\lambda)$ satisfies \eqref{G}, with $\Omega,\Gc,\P$ replaced respectively by $[0,1],\Bc([0,1]),\lambda$ (actually, every probability space $(E,\Ec,\Q)$, with $E$ uncountable, separable, complete metric space, $\Ec$ its Borel $\sigma$-algebra, $\Q$ an atomless probability, satisfies \eqref{G}, this follows e.g. from Corollary 7.16.1 in \cite{BertsekasShreve}).

Suppose now that the sub-$\sigma$-algebra $\Gc$ satisfies \eqref{G}. Denote by $\lambda$ the Lebesgue measure on $([0,1],\Bc([0,1]))$. Then, by the left-hand side equality in \eqref{G}, there exists a random variable $\Gamma^\Gc\colon(\Omega,\Gc)\rightarrow([0,1],\Bc(0,1))$ with distribution $\lambda$, that is with uniform distribution (so, in particular, $\Gamma^\Gc$ has an atomless distribution). On the other hand,  
given $\mu\in\mathscr P_{\text{\tiny$1$}}(\R)$ (resp. $\mu\in\mathscr P_{\text{\tiny$q$}}(\R^n)$) it is possible to find a random variable 
$\eta\colon[0,1]\rightarrow\R$ (resp. $\eta\colon[0,1]\rightarrow\R^n$) with distribution $\mu$. This implies that the random variable $\xi=\eta(\Gamma^\Gc)\colon\Omega\rightarrow\R$ (resp. $\xi=\eta(\Gamma^\Gc)\colon\Omega\rightarrow\R^n$) has also distribution $\mu$. 
As a matter of fact, we have:
\[
\P(\xi^{-1}(A)) \ = \ \P((\eta(\Gamma^\Gc))^{-1}(A)) \ = \ \P(\Gamma^\Gc\in\eta^{-1}(A)) \ = \ \lambda(\eta^{-1}(A)),
\]
for every $A\in\Bc(\R)$ (resp. $A\in\Bc(\R^n)$). Then, we see that the sub-$\sigma$-algebra $\bar\Gc:=\sigma(\Gamma^\Gc)\subset\Gc$ satisfies \eqref{G}. In other words, it is enough to replace $\Gc$ by the possibly smaller $\sigma$-algebra $\bar\Gc$.
\epR}
\end{Remark}

\begin{Remark}\label{R:ExH}
{\rm
Let $(E,\Ec)$ and $(H,\Hc)$ denote two Polish spaces endowed with their Borel $\sigma$-algebrae. Notice that the following result, stronger than \eqref{G}, holds true:
\begin{equation}\label{Kallenberg0}
\text{\emph{Given any probability $\pi$ on $(E,\Ec)$, there exists $\xi\colon(\Omega,\Gc)\rightarrow(E,\Ec)$ such that $\P_{\text{\tiny$\xi$}}=\pi$.}}
\end{equation}
Moreover, we have the following result (which will be used in the proof of Theorem \ref{prop:visco}):
\begin{align}\label{Kallenberg}
&\text{\emph{Given any probability $\pi$ on the product space $(E\times H,\Ec\otimes\Hc)$ and $\zeta\colon(\Omega,\Gc)\rightarrow(E,\Ec)$,}} \notag \\
&\text{\emph{with distribution $\P_{\text{\tiny$\zeta$}}$ equals to the marginal distribution of $\pi$ on $(E,\Ec)$,}} \\
&\text{\emph{there exists a random variable $\eta\colon(\Omega,\Gc)\rightarrow(H,\Hc)$ such that $\P_{\text{\tiny$(\zeta,\eta)$}}=\pi$.}} \notag
\end{align}
Proceeding along the same lines as in Remark \ref{R:G}, we see that statement \eqref{Kallenberg0} follows from Theorem 3.19 in \cite{Kallenberg} when $(G,\mathscr G)$ is $([0,1],\Bc([0,1]))$ and $\Gamma^\Gc$ has uniform distribution; in general, the claim follows simply recalling that all atomless Polish probability spaces are isomorphic (see for instance Corollary 7.16.1 in \cite{BertsekasShreve}).

Statement \eqref{Kallenberg} follows from Theorem 6.10 in \cite{Kallenberg} when $(G,\mathscr G)$ is $([0,1]\times[0,1],\Bc([0,1]\times[0,1]))$ and $\Gamma^\Gc$ has uniform distribution (proceeding as before, we deduce the result for the general case). As a matter of fact, given $\pi$ on $(E\times H,\Ec\otimes\Hc)$, by Theorem 3.19 in \cite{Kallenberg} there exists a random vector $(z,w)\colon([0,1]\times[0,1],\Bc([0,1]\times[0,1]))\rightarrow(E\times H,\Ec\otimes\Hc)$ such that $\P_{\text{\tiny$(z,w)$}}=\pi$. Now, since $\zeta$ is $\Gc$-measurable, by Doob's measurability theorem there exists a measurable map $\tilde z\colon[0,1]\times[0,1]\rightarrow E$ such that $\zeta=\tilde z(\Gamma^\Gc)$. Notice that $\tilde z$, as a random variable from $([0,1]\times[0,1],\Bc([0,1]\times[0,1]),\lambda\otimes\lambda)$ into $(E,\Ec)$, has the same distribution of $\zeta$, that is $\P_{\text{\tiny$\zeta$}}=\P_{\text{\tiny$\tilde z$}}$, which in turn coincides with the marginal distribution of $\pi$ on $(E,\Ec)$. So, in particular, $\P_{\text{\tiny$z$}}=\P_{\text{\tiny$\tilde z$}}$. We can now apply Theorem 6.10 in \cite{Kallenberg}, from which it follows the existence of $\tilde w\colon([0,1]\times[0,1],\Bc([0,1]\times[0,1]))\rightarrow(H,\Hc)$ such that $\P_{\text{\tiny$(\tilde z,\tilde w)$}}=\pi$. Define $\eta:=\tilde w(\Gamma^\Gc)$. Then, 
$\eta$ is a measurable map from $(\Omega,\Gc)$ into $(H,\Hc)$, moreover $\P_{\text{\tiny$(\zeta,\eta)$}}=\pi$, hence \eqref{Kallenberg} holds.}
\epR
\end{Remark}

Let $T>0$ be a finite time horizon and let $\Ac$ (resp. $\Bc$) be the family of admissible control processes for player I (resp. II), that is the set of all 
$(\Fc_s\vee\Gc)_s$-progressively measurable processes $\alpha\colon\Omega\times[0,T]\rightarrow A$ (resp. $\beta\colon\Omega\times[0,T]\rightarrow B$), where $A$ (resp. $B$) is a Polish space. We denote by $\rho_A$ (resp. $\rho_B$) a bounded metric on $A$ (resp. $B$) (notice that given a not necessarily bounded metric $d$ on a metric space $M$, the equivalent metric $d/(1+d)$ is bounded). Finally, we denote by $\mathscr P(A\times B)$ the family of all probability measures on $A\times B$, endowed with the topology of weak convergence.

The state equation of the McKean-Vlasov stochastic differential game is given by:
\begin{align}
X_s^{t,\xi,\alpha,\beta} \ &= \ \xi + \int_t^s \gamma\big(X_r^{t,\xi,\alpha,\beta},\P_{\text{\tiny$X_r^{t,\xi,\alpha,\beta}$}},\alpha_r,\beta_r,\P_{\text{\tiny$(\alpha_r,\beta_r)$}}\big)\,dr \label{State} \\
&\quad \ + \int_t^s \sigma\big(X_r^{t,\xi,\alpha,\beta},\P_{\text{\tiny$X_r^{t,\xi,\alpha,\beta}$}},\alpha_r,\beta_r,\P_{\text{\tiny$(\alpha_r,\beta_r)$}}\big)\,dW_r, \notag
\end{align}
for all $s\in[t,T]$, where $t\in[0,T]$, $\xi\in L^q(\Omega,\Fc_t\vee\Gc,\P;\R^n)$, $\alpha\in\Ac$, $\beta\in\Bc$. On the coefficients $b$ and $\sigma$, and also on the payoff functions $f$ and $g$ introduced below, we impose the following assumptions.

\vspace{2mm}

\noindent {\bf (A1)}
\begin{itemize}
\item [(i)] The maps $\gamma\colon\R^n\times\mathscr P_{\text{\tiny$q$}}(\R^n)\times A\times B\times\mathscr P(A\times B)\rightarrow\R^n$, 
$\sigma\colon\R^n\times\mathscr P_{\text{\tiny$q$}}(\R^n)\times A\times B\times\mathscr P(A\times B)\rightarrow\R^{n\times d}$, 
$f\colon\R^n\times\mathscr P_{\text{\tiny$q$}}(\R^n)\times A\times B\times\mathscr P(A\times B)\rightarrow\R$, $g\colon\R^n\times
\mathscr P_{\text{\tiny$q$}}(\R^n)\rightarrow\R$ 
are Borel measurable.
\item [(ii)] There exists a positive constant $L$  such that
\begin{align*}
|\gamma(x,\mu,a,b,\nu) - \gamma(x',\mu',a,b,\nu)| \ &\leq \ L\big(|x-x'| + \Wc_{\text{\tiny$q$}}(\mu,\mu')\big), \\
|\sigma(x,\mu,a,b,\nu)-\sigma(x',\mu',a,b,\nu)| \ &\leq \ L\big(|x-x'| + \Wc_{\text{\tiny$q$}}(\mu,\mu')\big), \\
|\gamma(0,\delta_0,a,b,\nu)| + |\sigma(0,\delta_0,a,b,\nu)| \ &\leq \ L, 
\\ |f(x,\mu,a,b,\nu)| + |g(x,\mu)| \ &\leq \ h(\|\mu\|_{\text{\tiny$q$}}) \big(1 + |x|^q\big),
\end{align*}
for all $(x,\mu),(x',\mu')\in\R^n\times\mathscr P_{\text{\tiny$q$}}(\R^n)$, $(a,b,\nu)\in A\times B\times\mathscr P(A\times B)$.
\end{itemize}

\begin{Remark}\label{R:JointLaw}
{\rm
In equation \eqref{State} we could also consider the case where $\gamma$ and $\sigma$ depend on the joint law $\P_{\text{\tiny$(X_r^{t,\xi,\alpha,\beta},\alpha_r,\beta_r)$}}$ rather than on the marginals $\P_{\text{\tiny$X_r^{t,\xi,\alpha,\beta}$}}$ and $\P_{\text{\tiny$(\alpha_r,\beta_r)$}}$. So, in particular, $\gamma=\gamma(x,a,b,\pi)$ and $\sigma=\sigma(x,a,b,\pi)$ for every $\pi\in\mathscr P_{\text{\tiny$q,0$}}(\R^n\times A\times B)$, the set of probability measures $\pi$ on the Borel $\sigma$-algebra of $\R^n\times A\times B$ with marginal $\pi_{_{|\R^n}}$ having finite $q$-th order moment. Notice however that such a generalization is only artificial, as a matter of fact under the Lipschitz assumption {\bf (A1)}-(ii), which now reads
\begin{align*}
|\gamma(x,a,b,\pi) - \gamma(x',a,b,\pi')| + |\sigma(x,a,b,\pi) - \sigma(x,a,b,\pi')| \ \leq \ L\big(|x-x'| + \Wc_{\text{\tiny$q$}}\big(\pi_{_{|\R^n}},\pi_{_{|\R^n}}'\big)\big), \\
\text{\emph{for all $(a,b)\in A\times B$ and $(x,\pi),(x',\pi')\in\R^n\times\mathscr P_{\text{\tiny$q,0$}}(\R^n\times A\times B)$, with $\pi_{_{|A\times B}}=\pi_{_{A\times B}}'$}},
\end{align*}
it follows that $\gamma(x,a,b,\pi)=\gamma(x,a,b,\pi')$ and $\sigma(x,a,b,\pi)=\sigma(x,a,b,\pi')$ whenever $\pi$ and $\pi'$ have the same marginals on $\R^n$ and $A\times B$. In other words, $\gamma=\gamma(x,a,b,\pi)$ and $\sigma=\sigma(x,a,b,\pi)$ depend only on the marginals of $\pi$ on $\R^n$ and $A\times B$.

Concerning the function $f$, we get to the same conclusion under the continuity assumption {\bf (A2)} stated below.
}
\epR
\end{Remark}

\begin{Lemma}
Under Assumption {\bf (A1)}, for any $t\in[0,T]$, $\xi\in L^q(\Omega,\Fc_t\vee\Gc,\P;\R^n)$, $\alpha\in\Ac$, $\beta\in\Bc$, there exists a unique (up to indistinguishability) continuous $(\Fc_s\vee\Gc)_s$-progressively measurable process $(X_s^{t,\xi,\alpha,\beta})_{s\in[t,T]}$ solution to \eqref{State}, satisfying
\begin{equation} \label{estimX}
\E\Big[\sup_{s\in[t,T]}\big|X_s^{t,\xi,\alpha,\beta}\big|^q\Big] \ \leq \ C_q\,\big(1 + \E[|\xi|^q]\big),
\end{equation}
for some positive constant $C_q$, independent of $t$, $\xi$, $\alpha$, $\beta$. Moreover, the flow property holds: for every $s\in[t,T]$,
\begin{equation}\label{flow}
X_r^{t,\xi,\alpha,\beta} \ = \ X_r^{s,X_s^{t,\xi,\alpha,\beta},\alpha,\beta}, \qquad \text{for all $r\in[s,T]$, $\P$-a.s.}
\end{equation}
and consequently
\begin{equation}\label{flowP}
\P_{\text{\tiny$X_r^{t,\xi,\alpha,\beta}$}} \ = \ \P_{\text{\tiny$X_r^{s,X_s^{t,\xi,\alpha,\beta},\alpha,\beta}$}}, \qquad \text{for all $r\in[s,T]$}.
\end{equation}
\end{Lemma}
\textbf{Proof.}
We report the proof only of \eqref{flow}-\reff{flowP}, the rest of the statement being standard. Notice that, by definition, the process $(X_r^{s,X_s^{t,\xi,\alpha,\beta},\alpha,\beta})_{r\in[s,T]}$ solves the following stochastic differential equation on $[s,T]$ with initial condition $X_s^{t,\xi,\alpha,\beta}$:
\begin{align*}
X_r \ &= \ X_s^{t,\xi,\alpha,\beta} + \int_s^r \gamma\big(X_z,\P_{\text{\tiny$X_z$}},\alpha_z,\beta_z,\P_{\text{\tiny$(\alpha_z,\beta_z)$}}\big)\,dz + \int_s^r \sigma\big(X_z,\P_{\text{\tiny$X_z$}},\alpha_z,\beta_z,\P_{\text{\tiny$(\alpha_z,\beta_z)$}}\big)\,dW_z,
\end{align*}
for all $r\in[s,T]$. On the other hand, recall from \eqref{State} that the process $X^{t,\xi,\alpha,\beta}$ solves the same equation on $[s,T]$, with identical initial condition at time $s$, that is $X_s^{t,\xi,\alpha,\beta}$. Hence, by pathwise uniqueness we conclude that $(X_r^{s,X_s^{t,\xi,\alpha,\beta},\alpha,\beta})_{r\in[s,T]}$ and $(X_r^{t,\xi,\alpha,\beta})_{r\in[s,T]}$ are indistinguishable, so that \eqref{flow} holds. We then deduce  the flow property 
\reff{flowP} on the probability law.  
\ep

\begin{Remark} \label{remlawxi}
{\rm  Notice that the (open-loop) control processes $\alpha$ $\in$ $\Ac$, $\beta$ $\in$ $\Bc$, are measurable with respect to $\Gc$, hence may depend on $\xi\in L^q(\Omega,\Fc_t\vee\Gc,\P;\R^n)$, and thus one cannot claim as in the uncontrolled case or when using feedback control that the law 
$\P_{\text{\tiny$X_s^{t,\xi,\alpha,\beta}$}}$  of $X_s^{t,\xi,\alpha,\beta}$, for $t\leq s\leq T$, depends on $\xi$ only through its distribution. 
}
\epR
\end{Remark}

\vspace{3mm}

The stochastic differential game has the following payoff functional:
\begin{equation}\label{J}
J(t,\xi,\alpha,\beta) \ = \ \E\bigg[\int_t^T f\big(X_r^{t,\xi,\alpha,\beta},\P_{\text{\tiny$X_r^{t,\xi,\alpha,\beta}$}},\alpha_r,\beta_r,\P_{\text{\tiny$(\alpha_r,\beta_r)$}}\big)\,ds + g\big(X_T^{t,\xi,\alpha,\beta},\P_{\text{\tiny$X_T^{t,\xi,\alpha,\beta}$}}\big)\bigg],
\end{equation}
for all $t\in[0,T]$, $\xi\in L^q(\Omega,\Fc_t\vee\Gc,\P;\R^n)$, $\alpha\in\Ac$, $\beta\in\Bc$. Notice that the payoff functional in \reff{J} is well-defined and finite by \reff{estimX} and the growth conditions on $f$ and $g$ in {\bf (A1)}(ii). We impose the following additional continuity conditions on the payoff functions $f$ and $g$.

\vspace{2mm}

\noindent {\bf (A2)} \quad The maps $g$ and $(x,\mu)$ $\in$  $\R^n\times\mathscr P_{\text{\tiny$q$}}(\R^n)$
$\mapsto$ $f(x,\mu,a,b,\nu)$ are  continuous, uniformly with respect to $(a,b,\nu)$ $\in$ $A\times B\times \mathscr P(A\times B)$, 
i.e. for any sequence $(x_m,\mu_m)_m$ in $\R^n\times\mathscr P_{\text{\tiny$q$}}(\R^n)$ converging to $(x,\mu)$ $\in$ 
$\R^n\times\mathscr P_{\text{\tiny$q$}}(\R^n)$, we have
\[
\sup_{_{(a,b,\nu) \in A\times B\times \mathscr P(A\times B)}} \big| f(x_m,\mu_m,a,b,\nu) - f(x,\mu,a,b,\nu) \big| \ + \ |g(x_m,\mu_m) - g(x,\mu)| \ \overset{m\rightarrow\infty}{\longrightarrow} \ 0. 
\]

\vspace{2mm}

We define the upper and lower value functions of the stochastic differential game as in Definition 1.4 of \cite{FlemingSouganidis}. In order to do it, we need to introduce the concept of (non-anticipative) strategy (see Definition 1.3 in \cite{FlemingSouganidis}).

\begin{Definition}\label{D:Str}
\quad
\begin{itemize}
\item A strategy $\alpha[\cdot]$ for player I  is a map $\alpha[\cdot]\colon\Bc\rightarrow\Ac$ satisfying the non-anticipativity property:
\[
\P\big(\beta_r\,=\,\beta_r',\text{ for a.e. }r\in[0,t]\big) \ = \ 1 \quad \Longrightarrow 
\quad \P\big(\alpha[\beta]_r\,=\,\alpha[\beta']_r,\text{ for a.e. }r\in[0,t]\big) \ = \ 1,
\]
for every $t\in[0,T]$ and any $\beta,\beta'\in\Bc$. We denote by $\Ac_{\textup{\tiny str}}$ the family of all strategies for player I.
\item A strategy $\beta[\cdot]$ for player II  is a map $\beta[\cdot]\colon\Ac\rightarrow\Bc$ satisfying the non-anticipativity property:
\[
\P\big(\alpha_r\,=\,\alpha_r',\text{ for a.e. }r\in[0,t]\big) \ = \ 1 \quad \Longrightarrow 
\quad \P\big(\beta[\alpha]_r\,=\,\beta[\alpha']_r,\text{ for a.e. }r\in[0,t]\big) \ = \ 1,
\]
for every $t\in[0,T]$ and any $\alpha,\alpha'\in\Ac$. We denote by $\Bc_{\textup{\tiny str}}$ the family of all strategies for player II.
\end{itemize}
\end{Definition}

The lower value function of the stochastic differential game (SDG) is given by
\[
v(t,\xi) \ = \ \inf_{\beta[\cdot]\in\Bc_{\textup{\tiny str}}} \sup_{\alpha\in\Ac} J(t,\xi,\alpha,\beta[\alpha]), \qquad \text{for all }t\in[0,T],\,\xi\in L^q(\Omega,\Fc_t\vee\Gc,\P;\R^n).
\]
On the other hand, the upper value function of the stochastic differential game is given by
\[
u(t,\xi) \ = \ \sup_{\alpha[\cdot]\in\Ac_{\textup{\tiny str}}} \inf_{\beta\in\Bc} J(t,\xi,\alpha[\beta],\beta), \qquad \text{for all }t\in[0,T],\,\xi\in L^q(\Omega,\Fc_t\vee\Gc,\P;\R^n).
\]
By estimate \reff{estimX} and the growth conditions in {\bf (A1)}(ii) on $f$ and $g$, we easily see that the value functions $v$ and $u$ satisfy the growth condition
\beqs
|v(t,\xi)| + |u(t,\xi)| & \leq & C_q\,h\big(C_q(1+\|\mu\|_{_q})\big)\,\big(1 + \E|\xi|^q\big), \qquad  t\in[0,T],\,\xi\in L^q(\Omega,\Fc_t\vee\Gc,\P;\R^n),
\enqs
with $C_q$ as in estimate \eqref{estimX}, and $\mu$ $=$ $\P_{_\xi}$.

\section{Properties of the value functions} \label{secprop}

The main goal of this section is to prove that the lower and upper value functions   $v$ and $u$ are law-invariant, i.e.,  
depend on $\xi$ only via its distribution. We mention that this is a non-trivial issue, as 
such property does not hold in general  for the payoff functional $J$, see Remark \ref{remlawxi}. 

\begin{Proposition}\label{P:Lift}
Under Assumption {\bf (A1)}, for every $t\in[0,T]$ we have
\[
v(t,\xi) \ = \ v(t,\tilde\xi), \qquad\qquad u(t,\xi) \ = \ u(t,\tilde\xi),
\]
for any $\xi,\tilde\xi\in L^q(\Omega,\Fc_t\vee\Gc,\P;\R^n)$, with $\P_{\text{\tiny$\xi$}}=\P_{\text{\tiny$\tilde\xi$}}$.
\end{Proposition}
\textbf{Proof.}
We prove the result only for the lower value function, as the proof for the upper value function can be done proceeding along the same lines.

Fix $t\in[0,T]$, $\xi,\tilde\xi\in L^q(\Omega,\Fc_t\vee\Gc,\P;\R^n)$, with $\mu:=\P_{\text{\tiny$\xi$}}=\P_{\text{\tiny$\tilde\xi$}}$. Our aim is to prove that:
\begin{align}\label{JalphaJtildealpha}
&\text{\emph{Given any $\alpha\in\Ac$ and $\tilde\beta[\cdot]\in\Bc_{\textup{\tiny str}}$, there exist $\tilde\alpha\in\Ac$ and $\beta[\cdot]\in\Bc_{\textup{\tiny str}}$, with $\tilde\alpha$ (resp. $\beta[\cdot]$)}} \notag \\
&\text{\emph{independent of $\tilde\beta[\cdot]$ (resp. $\alpha$), depending only on $\xi$, $\tilde\xi$, $\alpha$ (resp. $\xi$, $\tilde\xi$, $\tilde\beta[\cdot]$), such that:}} \\
&(\xi,\alpha,\beta[\alpha],W)\text{\emph{ has the same law as }}(\tilde\xi,\tilde\alpha,\tilde\beta[\tilde\alpha],W)\text{\emph{, so that }}J(t,\xi,\alpha,\beta[\alpha]) \ = \ J(t,\tilde\xi,\tilde\alpha,\tilde\beta[\tilde\alpha]). \notag
\end{align}
Notice that statement \eqref{JalphaJtildealpha} is equivalent to the existence of two maps $\Psi_1\colon\Ac\rightarrow\Ac$ and $\Psi_2\colon\Bc_{\textup{\tiny str}}\rightarrow\Bc_{\textup{\tiny str}}$ 
such that for any $\alpha\in\Ac$ and $\tilde\beta[\cdot]\in\Bc_{\textup{\tiny str}}$, the quadruple $(\xi,\alpha,$ $\Psi_2(\tilde\beta[\cdot])[\alpha],W)$ has the same law as $(\tilde\xi,\Psi_1(\alpha),\tilde\beta[\Psi_1(\alpha)],W)$, so that $J(t,\xi,\alpha,\Psi_2(\tilde\beta[\cdot])[\alpha])$ $=$ $J(t,\tilde\xi,\Psi_1(\alpha),\tilde\beta[\Psi_1(\alpha)])$.

Observe that the claim follows if \eqref{JalphaJtildealpha} holds true. Indeed, for any fixed $\tilde\beta[\cdot]\in\Bc_{\textup{\tiny str}}$ we have
\[
\sup_{\alpha\in\Ac} J(t,\xi,\alpha,\Psi_2(\tilde\beta[\cdot])[\alpha]) \ = \ \sup_{\alpha\in\Ac} J(t,\tilde\xi,\Psi_1(\alpha),\tilde\beta[\Psi_1(\alpha)]) \ \leq \ \sup_{\tilde\alpha\in\Ac} J(t,\tilde\xi,\tilde\alpha,\tilde\beta[\tilde\alpha]).
\]
Taking the infimum over $\tilde\beta[\cdot]$ in $\Bc_{\textup{\tiny str}}$, we obtain
\begin{align*}
v(t,\xi) \ = \ \inf_{\beta[\cdot]\in\Bc_{\textup{\tiny str}}} \sup_{\alpha\in\Ac} J(t,\xi,\alpha,\beta[\alpha]) \ &\leq \ \inf_{\tilde\beta[\cdot]\in\Bc_{\textup{\tiny str}}} \sup_{\alpha\in\Ac} J(t,\xi,\alpha,\Psi_2(\tilde\beta[\cdot])[\alpha]) \\
&\leq \ \inf_{\tilde\beta[\cdot]\in\Bc_{\textup{\tiny str}}} \sup_{\tilde\alpha\in\Ac} J(t,\tilde\xi,\tilde\alpha,\tilde\beta[\tilde\alpha]) \ = \ v(t,\tilde\xi).
\end{align*}
Interchanging the roles of $\xi$ and $\tilde\xi$, we get the other inequality $v(t,\xi)\geq v(t,\tilde\xi)$, from which the claim follows.

It remains to prove statement \eqref{JalphaJtildealpha}. We split its proof into four steps. Notice that \textbf{Step II} is similar to the proof of Proposition 2.2 in \cite{BBK}, while \textbf{Step IV} can be alternatively addressed using techniques from the proof of Proposition 3.4 in \cite{BBK} (see \textbf{Step IV} below for more details).

\vspace{2mm}

\noindent\textbf{Step I.} \emph{Reduction to a canonical setting.} Denote by $E_t^W:=C_0([0,t];\R^n)$ the set of $\R^n$-valued continuous paths on $[0,t]$ starting at the origin at time $0$. We endow $E_t^W$ with the uniform topology, so that $E_t^W$ becomes a Polish space (we denote its Borel $\sigma$-algebra by $\Ec_t^W$). We also denote by $\P_t^W$ the Wiener measure on $(E_t^W,\Ec_t^W)$ (recall that the Wiener measure is atomless). Now, consider the filtration $\F^o=(\Fc_s^o)_{s\geq0}$ generated by the Brownian motion $W$. Notice that there exists a random variable $\Gamma_t^W\colon(\Omega,\Fc_t^o)\rightarrow(E_t^W,\Ec_t^W)$ with distribution $\P_t^W$ and such that $\Fc_t^o=\sigma(\Gamma_t^W)$. On the other hand, we recall that, by assumption, there exists a random variable $\Gamma^\Gc\colon(\Omega,\Gc)\rightarrow(G,\mathscr G)$, taking values in some Polish space $G$ with Borel $\sigma$-algebra $\mathscr G$, such that $\Gamma^\Gc$ has an atomless distribution and $\Gc=\sigma(\Gamma^\Gc)$. Hence, we deduce that there exists a random variable $\Gamma_t\colon(\Omega,\Fc_t^o\vee\Gc)\rightarrow(E,\Ec)$, taking values in some Polish space $E$ with Borel $\sigma$-algebra $\Ec$, such that $\Gamma_t$ has an atomless distribution and $\Fc_t^o\vee\Gc=\sigma(\Gamma_t)$. Finally, recalling that all atomless Polish probability spaces are isomorphic, we can suppose that the probability space $(E,\Ec,\P_{\text{\tiny$\Gamma_t$}})$, where $\P_{\text{\tiny$\Gamma_t$}}$ denotes the distribution of $\Gamma_t$, is given by the space $([0,1],\Bc([0,1]),\lambda)$, where $\lambda$ is the Lebesgue measure on $[0,1]$. So, in particular, $\Gamma_t\colon\Omega\rightarrow[0,1]$ and has uniform distribution.

\vspace{2mm}

\noindent\textbf{Step II.} \emph{Canonical representation of $\xi$ and $\alpha$.} Fix $\alpha\in\Ac$ and $\tilde\beta[\cdot]\in\Bc_{\textup{\tiny str}}$.

\vspace{1mm}

\noindent\emph{Representation of $\xi$.} Since $\xi$ is $\Fc_t\vee\Gc$-measurable, by Doob's measurability theorem it follows that
\[
\xi \ = \ \chi(\Gamma_t), \qquad \P\text{-a.s.}
\]
for some measurable function $\chi\colon([0,1],\Bc([0,1]))\rightarrow(\R^n,\Bc(\R^n))$. The equality $\xi=\chi(\Gamma_t)$ holds $\P$-a.s. since $\Fc_t=\Fc_t^o\vee\Nc$. Notice that we can suppose $\chi$ to be surjective. As a matter of fact, if this is not the case, it is enough to modify $\chi$ on the set $\mathscr C\backslash\{0,1\}$ (where $\mathscr C$ is the Cantor set), replacing $\chi$ for instance by the composition of the Cantor function from $\mathscr C\backslash\{0,1\}$ to $(0,1)$ with a continuous map from $(0,1)$ to $\R^n$. The $\chi$ so constructed remains a Borel measurable function. Moreover, we still have
\[
\xi \ = \ \chi(\Gamma_t), \qquad \P\text{-a.s.}
\]

\vspace{1mm}

\noindent\emph{Representation of $\alpha$.} Similarly, the map $\alpha\colon\Omega\times[0,T]\rightarrow A$ is $Prog(\F^t)\vee((\Fc_t\vee\Gc)\otimes\{\emptyset,[0,T]\})$-measurable, where 
$Prog(\F^t)$ denotes the progressive $\sigma$-algebra on $\Omega\times[0,T]$ relative to the filtration $\F^t$, while $\{\emptyset,[0,T]\}$ is the trivial $\sigma$-algebra on $[0,T]$. Then, by a slight generalization of Doob's measurability theorem (which can be proved using the monotone class theorem), it follows that $\alpha$ has the form 
$\alpha_s=a_s(\cdot,\Gamma_t(\cdot))$, $\forall\,s\in[0,T]$, $\P$-a.s., for some $Prog(\F^t)\otimes\Bc([0,1])$-measurable function $a=a_s(\omega,y)\colon\Omega\times[0,T]\times[0,1]\rightarrow A$. As before, the fact that the equality $\alpha_s=a_s(\cdot,\Gamma_t(\cdot))$, $\forall\,s\in[0,T]$, holds $\P$-a.s. (so, in particular, 
$(\alpha_s)_s$ and $(a_s(\cdot,\Gamma_t(\cdot)))_s$ are $\P$-indistinguishable) follows from the fact that $\Fc_t=\Fc_t^o\vee\Nc$.

\vspace{2mm}

\noindent\textbf{Step III.} \emph{The random variable $\tilde\Gamma_t$.} Notice that \eqref{JalphaJtildealpha} follows if we prove the following:
\begin{align}\label{tildeGamma}
&\text{\emph{$\exists$ a random variable $\tilde\Gamma_t\colon(\Omega,\Fc_t\vee\Gc)\rightarrow([0,1],\Bc([0,1]))$ such that:}} \\
&\text{\emph{$\tilde\Gamma_t$ has the same distribution of $\Gamma_t$, $\sigma(\tilde\Gamma_t)\vee\Nc=\sigma(\Gamma_t)\vee\Nc=\Fc_t\vee\Gc$, and }} \tilde\xi \ = \ \chi(\tilde\Gamma_t),\,\P\text{\emph{-a.s.}} \notag
\end{align}
Observe that $\tilde\Gamma_t$ is allowed to be $\Fc_t\vee\Gc$-measurable (not necessarily just $\Fc_t^o\vee\Gc$-measurable). Suppose that \eqref{tildeGamma} holds. Since $\sigma(\tilde\Gamma_t)\vee\Nc=\Fc_t\vee\Gc$, we can find a canonical representation of $\tilde\beta[\cdot]$ in terms of $\tilde\Gamma_t$. More precisely, since for every $\tilde\alpha\in\Ac$ we have that $\tilde\beta[\tilde\alpha]$ is an element of $\Bc$, proceeding as in \textbf{Step II} for the proof of the canonical representation of $\alpha$, we deduce that $\tilde\beta[\tilde\alpha]$ has the form $\tilde\beta[\tilde\alpha]_s=\tilde b_s^{\tilde\alpha}(\cdot,\tilde\Gamma_t(\cdot))$, $\forall\,s\in[0,T]$, $\P$-a.s., for some $Prog(\F^t)\otimes\Bc([0,1])$-measurable function $\tilde b^{\tilde\alpha}=\tilde b_s^{\tilde\alpha}(\omega,y)\colon\Omega\times[0,T]\times[0,1]\rightarrow B$. Now, define
\[
\big(\Psi_1(\alpha) \, =\big) \ \tilde\alpha \ := \ a_\cdot(\cdot,\tilde\Gamma_t(\cdot)), \qquad\qquad \big(\Psi_2(\tilde\beta[\cdot])[\cdot] \, =\big) \ \beta[\cdot] \ := \ \tilde b_\cdot^\cdot(\cdot,\Gamma_t(\cdot)).
\]
Notice that $\tilde\alpha\in\Ac$ and $\beta[\cdot]\in\Bc_{\textup{\tiny str}}$. We also notice that $(\xi,\alpha,\beta[\alpha],W)$ has the same law as $(\tilde\xi,\tilde\alpha,\tilde\beta[\tilde\alpha],W)$. So, in particular, $J(t,\xi,\alpha,\beta[\alpha])=J(t,\tilde\xi,\tilde\alpha,\tilde\beta[\tilde\alpha])$, that is \eqref{JalphaJtildealpha} holds.

\vspace{2mm}

\noindent\textbf{Step IV.} \emph{Proof of \eqref{tildeGamma}.} By the Jankov-von Neumann measurable selection theorem (see for instance Theorem 18.22 and, in particular, Corollary 18.23 in \cite{Selection}), it follows that $\chi$ admits an analytically measurable right-inverse, denoted by $\zeta\colon\R^n\rightarrow [0,1]$, which satisfies:
\begin{itemize}
\item[1)] $\chi(\zeta(y))=y$, for any $y\in\R^n$;
\item[2)] $\chi^{-1}(\zeta^{-1}(\mathscr B))=\mathscr B$, for any subset $\mathscr B$ of $[0,1]$;
\item[3)] $\zeta^{-1}(\mathscr B)$ is analytically measurable in $\R^n$ for each Borel subset $\mathscr B$ of $[0,1]$. Recalling that every analytic subset of $\R^n$ is universally measurable (see e.g. Theorem 12.41 in \cite{Selection}), it follows that $\zeta^{-1}(\mathscr B)\in\mathscr L(\R^n)$, the Lebesgue $\sigma$-algebra on $\R^n$. Hence $\zeta$ is a measurable function from $(\R^n,\mathscr L(\R^n))$ into $([0,1],\Bc([0,1]))$.
\end{itemize}
Now, define  
\[
\tilde\Gamma_t \ := \ \zeta(\tilde\xi),
\]
and let us prove that $\tilde\Gamma_t$ satisfies \eqref{tildeGamma} (notice that in the proof of Proposition 3.4 in \cite{BBK}, $\tilde\Gamma_t$ is constructed in a different and more direct way, namely by means of the regular conditional cumulated distribution of $\Gamma_t$ given both $\xi$ and $(W_s)_{0\leq s\leq t}$, see \cite{BBK} for all the details). 

We begin noting that, since $\zeta$ is $\Bc([0,1])/\mathscr L(\R^n)$-measurable, and also the $\sigma$-algebra $\Fc_t\vee\Gc$ is $\P$-complete, it follows that $\tilde\Gamma_t$ is a measurable function from $(\Omega,\Fc_t\vee\Gc)$ into $([0,1],\Bc([0,1]))$. Let us now prove that $\tilde\Gamma_t$ has the same distribution of $\Gamma_t$. Fix a Borel subset $\mathscr B$ of $[0,1]$. Then
\[
\P(\tilde\Gamma_t\in\mathscr B) \ = \ \P(\zeta(\tilde\xi)\in\mathscr B) \ = \ \P(\tilde\xi\in\zeta^{-1}(\mathscr B)).
\]
Recalling that $\tilde\xi$ has the same distribution of $\xi$, and also that $\xi=\chi(\Gamma_t)$, we obtain
\begin{align*}
\P(\tilde\xi\in\zeta^{-1}(\mathscr B)) \ = \ \P(\chi(\Gamma_t)\in\zeta^{-1}(\mathscr B)) \ = \ \P(\Gamma_t\in\chi^{-1}(\zeta^{-1}(\mathscr B))).
\end{align*}
By item 2), we know that $\chi^{-1}(\zeta^{-1}(\mathscr B))=\mathscr B$, hence
\[
\P(\tilde\Gamma_t\in\mathscr B) \ = \ \P(\Gamma_t\in\mathscr B). 
\]
This proves that $\tilde\Gamma_t$ has the same distribution of $\Gamma_t$. Moreover, by item 1) we have $\chi(\tilde\Gamma_t)=\chi(\zeta(\tilde\xi))=\tilde\xi$. 
It remains to prove the equality $\sigma(\tilde\Gamma_t)\vee\Nc=\sigma(\Gamma_t)\vee\Nc$.

Similarly to $\xi$, since $\tilde\xi$ is $\Fc_t\vee\Gc$-measurable, by Doob's measurability theorem there exists a measurable function $\tilde\chi\colon([0,1],\Bc([0,1]))\rightarrow(\R^n,\Bc(\R^n))$ such that $\tilde\xi=\tilde\chi(\Gamma_t)$, $\P$-a.s.. Hence
\[
\tilde\Gamma_t \ = \ \zeta(\tilde\chi(\Gamma_t)), \qquad \P\text{-a.s.}
\]
So, in particular, $\sigma(\tilde\Gamma_t)\subset\sigma(\Gamma_t)\vee\Nc$. It remains to prove that $\sigma(\Gamma_t)\subset\sigma(\tilde\Gamma_t)\vee\Nc$. Notice that $\zeta\circ\tilde\chi$ is a measurable function from $([0,1],\mathscr L([0,1]))$ into $([0,1],\Bc([0,1]))$. Then, it is well-known that there exists a Borel measurable function $\phi\colon[0,1]\rightarrow[0,1]$ such that $\zeta\circ\tilde\chi=\phi$, $\lambda$-a.e. (we need to consider a Borel measurable version of $\zeta\circ\tilde\chi$ in order to use Lemma \ref{L:AppTechn}, which in turn relies on the Jankov-von Neumann measurable selection theorem). Hence
\[
\tilde\Gamma_t \ = \ \phi(\Gamma_t), \qquad \P\text{-a.s.}
\]
Then, by Lemma \ref{L:AppTechn} it follows that there exists a Borel measurable function $\rho\colon[0,1]\rightarrow[0,1]$ such that $\rho(\phi)(y)=y$, $\lambda$-a.e., so that
\[
\Gamma_t \ = \ \rho(\tilde\Gamma_t), \qquad \P\text{-a.s.},
\]
from which we deduce the inclusion $\sigma(\Gamma_t)\subset\sigma(\tilde\Gamma_t)\vee\Nc$. This concludes the proof.
\ep

\begin{Definition} \label{definverselift}
By Proposition \ref{P:Lift}, we define the inverse-lifted functions of $v$ and $u$, respectively:
\[
\mathscr V(t,\mu) \ = \ v(t,\xi), \qquad\qquad \mathscr U(t,\mu) \ = \ u(t,\xi), \qquad \text{for every }(t,\mu)\in[0,T]\times\mathscr P_{\text{\tiny$q$}}(\R^n),
\]
for any $\xi\in L^q(\Omega,\Fc_t\vee\Gc,\P;\R^n)$, with $\P_{\text{\tiny$\xi$}}=\mu$.
\end{Definition}

We end this section proving the continuity of the value functions.

\begin{Proposition}\label{P:Cont}
Let Assumptions {\bf (A1)}-{\bf (A2)} hold. The function $(t,\xi)$ $\mapsto$ $J(t,\xi,\alpha,\beta)$ is continuous on $\Dc$ $:=$ 
$\{ t \in [0,T], \xi \in  L^q(\Omega,\Fc_t\vee\Gc,\P;\R^n)\}$, uniformly  with respect to $(\alpha,\beta)$ $\in$ $\Ac\times\Bc$, and consequently the value functions $v$, $u$ are continuous on $\Dc$. 
\end{Proposition}
\textbf{Proof.} (1)  Fix $0\leq t\leq s \leq T$, $\xi,\zeta$ $\in$ $L^q(\Omega,\Fc_t\vee\Gc,\P;\R^n)$, $\alpha$ $\in$ $\Ac$, $\beta$ $\in$ $\Bc$. 
By definition of Wasserstein distance, we have
\beq \label{estimWasser}
\sup_{s \leq r \leq T} \Wc_{_q}^q\big( \P_{_{X_r^{t,\xi,\alpha,\beta}}},  \P_{_{X_r^{s,\zeta,\alpha,\beta}}} \big) & \leq &  
\E \Big[  \sup_{s \leq r \leq T} \big|X_r^{t,\xi,\alpha,\beta} -  X_r^{s,\zeta,\alpha,\beta} \big|^q \Big]. 
\enq
From the state equation \reff{State}, and using standard arguments involving Burkholder-Davis-Gundy inequalities, \reff{estimWasser}, and Gronwall lemma, under the Lipschitz condition in {\bf(A1)}(ii), we obtain the following estimate similar to the ones  for controlled diffusion processes (see Theorem 5.9 and Corollary 5.10, Chapter 2, in \cite{Krylov}):   
\beq \label{estimX1}
\E\Big[\sup_{s \leq  r \leq T} \big|X_r^{t, \xi, \alpha,\beta}-X_r^{s, \zeta, \alpha,\beta}\big|^q \Big] &  \leq &  C \big(\E|\xi-\zeta|^q +(1 + \E|\xi|^q + \E|\zeta|^q)|s-t| \big),
\enq
for some constant $C$ (independent of $t,s,\xi,\zeta, \alpha,\beta$).

\noindent (2)  
Fix $t \in [0,T]$ and $\xi\in L^q(\Omega,\Fc_t\vee\Gc,\P;\R^n)$. Consider a sequence $(t_m)_m$ $\subset$ 
$[0,T]$, $(\xi_m)_m\subset L^q(\Omega,\Fc_{t_m}\vee\Gc,\P;\R^n)$ such that
$t_m$ $\rightarrow$ $t$, and $\xi_m$ $\rightarrow$ $\xi$ in $L^q$ as $m$ goes to infinity.  We then have  for all $\alpha$ $\in$ $\Ac$, $\beta$ $\in$ $\Bc$, 
\begin{align*}
&\big| J(t_m,\xi_m,\alpha,\beta) - J(t,\xi,\alpha,\beta) \big| \\
&\leq \ \E \bigg[ \int_t^{t_m}  \big|f\big(X_r^{t,\xi,\alpha,\beta},\P_{\text{\tiny$X_r^{t,\xi,\alpha,\beta}$}},\alpha_r,\beta_r,\P_{\text{\tiny$(\alpha_r,\beta_r)$}}\big) \big|  dr \\
&\quad \ + \int_t^T   \big|f\big(X_r^{t,\xi,\alpha,\beta},\P_{\text{\tiny$X_r^{t,\xi,\alpha,\beta}$}},\alpha_r,\beta_r,\P_{\text{\tiny$(\alpha_r,\beta_r)$}}\big) 
-  f\big(X_r^{t_m,\xi_m,\alpha,\beta},\P_{\text{\tiny$X_r^{t_m,\xi_m,\alpha,\beta}$}},\alpha_r,\beta_r,\P_{\text{\tiny$(\alpha_r,\beta_r)$}}\big) \big|  dr \\
&\quad \ + \big| g\big(X_T^{t,\xi,\alpha,\beta},\P_{\text{\tiny$X_T^{t,\xi,\alpha,\beta}$}}\big) - g\big(X_T^{t_m,\xi_m,\alpha,\beta},\P_{\text{\tiny$X_T^{t_m,\xi_m,\alpha,\beta}$}}\big) \big| \bigg] \\
&\leq \ C_q\,h\big(C_q(1+\E|\xi|^q)\big)\,\big(1 + \E|\xi|^q\big)\,|t_m-t|  \\
&\quad \ + \E\bigg[\!\int_t^T \!\!\!\!\sup_{_{(r,a,b,\nu) \in [0,T]\times A\times B\times \mathscr P(A\times B)}} 
\!\!\!\!\!\!\!\!\!\!\!\!\!\!\!\big|f\big(X_r^{t,\xi,\alpha,\beta},\P_{\text{\tiny$X_r^{t,\xi,\alpha,\beta}$}},a,b,\nu\big)\!-\!f\big(X_r^{t_m,\xi_m,\alpha,\beta},\P_{\text{\tiny$X_r^{t_m,\xi_m,\alpha,\beta}$}},a,b,\nu\big) \big|   dr \!\bigg]  \\
&\quad \ + \E\Big[ \big|g\big(X_T^{t,\xi,\alpha,\beta},\P_{\text{\tiny$X_T^{t,\xi,\alpha,\beta}$}}\big) - g\big(X_T^{t_m,\xi_m,\alpha,\beta},\P_{\text{\tiny$X_T^{t_m,\xi_m,\alpha,\beta}$}}\big) \big| \Big].
\end{align*}
From \reff{estimWasser}-\reff{estimX1}, and the continuity assumption {\bf (A2)} on $f$ and $g$, we deduce that 
\beqs
\sup_{\alpha\in\Ac,\beta\in\Bc}|J(t_m,\xi_m,\alpha,\beta) - J(t,\xi,\alpha,\beta)| \ \underset{m\rightarrow\infty}{\longrightarrow} \ 0,
\enqs
which implies 
\[
v(t_m,\xi_m) \ \underset{m\rightarrow\infty}{\longrightarrow} \ v(t,\xi), \qquad\qquad u(t_m,\xi_m) \ \underset{m\rightarrow\infty}{\longrightarrow} \ u(t,\xi),
\]
from which the claim follows.
\ep

\begin{Corollary}
Let Assumptions {\bf (A1)}-{\bf (A2)} hold. The inverse-lifted functions $\mathscr V$ and $\mathscr U$ are continuous on $[0,T]\times \mathscr P_{\text{\tiny$q$}}(\R^n)$.
\end{Corollary}
\textbf{Proof.} The claim follows directly from the continuity of the value functions $v$ and $u$ in Proposition \ref{P:Cont}, and also by Skorohod's representation theorem on the Wasserstein space (see Lemma A.1 in \cite{bayetal}).
\ep

\section{Dynamic Programming Principle} \label{secDPP}

The main result of this section is the statement and proof of the dynamic programming principle (DPP) for the lower and upper value functions of the two-player zero-sum McKean-Vlasov stochastic differential game.

\begin{Theorem}\label{T:DPP}
Under Assumption {\bf (A1)}, we have
\begin{align}\label{DPPv}
v(t,\xi) \ = \ \inf_{\beta[\cdot]\in\Bc_{\textup{\tiny str}}} \sup_{\alpha\in\Ac} \bigg\{ \E\Big[\int_t^s f\big(X_r^{t,\xi,\alpha,\beta[\alpha]},\P_{\text{\tiny$X_r^{t,\xi,\alpha,\beta[\alpha]}$}},\alpha_r,\beta[\alpha]_r,\P_{\text{\tiny$(\alpha_r,\beta[\alpha]_r)$}}\big)\,dr&\Big] \\
\;\;\; +  \;\;  v\big(s,X_s^{t,\xi,\alpha,\beta[\alpha]}\big)&\bigg\} \notag
\end{align}
and
\begin{align}\label{DPPu}
u(t,\xi) \ = \ \sup_{\alpha[\cdot]\in\Ac_{\textup{\tiny str}}} \inf_{\beta\in\Bc} \bigg\{ \E\Big[\int_t^s f\big(X_r^{t,\xi,\alpha[\beta],\beta},\P_{\text{\tiny$X_r^{t,\xi,\alpha[\beta],\beta}$}},\alpha[\beta]_r,\beta_r,\P_{\text{\tiny$(\alpha[\beta]_r,\beta_r)$}}\big)\,dr&\Big] \\
+ \, u\big(s,X_s^{t,\xi,\alpha[\beta],\beta}\big)&\bigg\}, \notag
\end{align}
for all $t,s\in[0,T]$, with $t\leq s$, and for every $\xi\in L^q(\Omega,\Fc_t\vee\Gc,\P;\R^n)$.
\end{Theorem}
\textbf{Proof.}
We prove the dynamic programming principle \eqref{DPPv} for the lower value function $v$, the proof of \eqref{DPPu} being similar.

For any $t\in[0,T]$, $\xi\in L^q(\Omega,\Fc_t\vee\Gc,\P;\R^n)$, $\eta\in L^1(\Omega,\Fc_t\vee\Gc,\P;\R)$, $\alpha\in\Ac$, $\beta\in\Bc$, consider the stochastic process $(\tilde X_s^{t,\xi,\eta,\alpha,\beta})_{s\in[t,T]}$ defined as
\[
\tilde X_s^{t,\xi,\eta,\alpha,\beta} \ := \ \eta + \int_t^s f\big(X_r^{t,\xi,\alpha,\beta},\P_{\text{\tiny$X_r^{t,\xi,\alpha,\beta}$}},\alpha_r,\beta_r,\P_{\text{\tiny$(\alpha_r,\beta_r)$}}\big)\,dr, 
\;\;\;  t \leq s \leq T. 
\]
Notice that, from identities \eqref{flow}-\eqref{flowP}, we deduce the following flow property: for every $s\in[t,T]$,
\begin{equation}\label{flow_tildeX}
\tilde X_r^{t,\xi,\eta,\alpha,\beta} \ = \ \tilde X_r^{s,X_s^{t,\xi,\alpha,\beta},\tilde X_s^{t,\xi,\eta,\alpha,\beta},\alpha,\beta}, \qquad \text{for all $r\in[s,T]$, $\P$-a.s.}
\end{equation}
Now, we observe that
\begin{align*}
J(t,\xi,\alpha,\beta) &= \ \E\big[\tilde X_T^{t,\xi,0,\alpha,\beta} + g\big(X_T^{t,\xi,\alpha,\beta},\P_{\text{\tiny$X_T^{t,\xi,\alpha,\beta}$}}\big)\big] \ = \ G\big(X_T^{t,\xi,\alpha,\beta},\tilde X_T^{t,\xi,0,\alpha,\beta}\big),
\end{align*}
where $G\colon L^q(\Omega,\Fc_T\vee\Gc,\P;\R^n)\times L^1(\Omega,\Fc_T\vee\Gc,\P;\R)\rightarrow\R$ is defined as
\[
G(\xi,\eta) \ := \ \E\big[\eta + g\big(\xi,\P_{\text{\tiny$\xi$}}\big)\big], \qquad \forall\,(\xi,\eta)\in L^q(\Omega,\Fc_T\vee\Gc,\P;\R^n)\times L^1(\Omega,\Fc_T\vee\Gc,\P;\R).
\]
Then, the lower value function of the stochastic differential game is given by
\[
v(t,\xi) \ = \ \inf_{\beta[\cdot]\in\Bc_{\textup{\tiny str}}} \sup_{\alpha\in\Ac} G\big(X_T^{t,\xi,\alpha,\beta[\alpha]},\tilde X_T^{t,\xi,0,\alpha,\beta[\alpha]}\big).
\]
Let
\[
V(t,\xi,\eta) \ = \ \inf_{\beta[\cdot]\in\Bc_{\textup{\tiny str}}} \sup_{\alpha\in\Ac} G\big(X_T^{t,\xi,\alpha,\beta[\alpha]},\tilde X_T^{t,\xi,\eta,\alpha,\beta[\alpha]}\big).
\]
Notice that the following relation holds between $v$ and $V$:
\begin{equation}\label{tildeV=V}
V(t,\xi,\eta) \ = \ v(t,\xi) + \E[\eta].
\end{equation}
By \eqref{tildeV=V}, we see that the dynamic programming principle \eqref{DPPv} for $v$ holds if and only if the following dynamic programming principle for $V$ holds:
\begin{equation}\label{DPPtildeV}
V(t,\xi,\eta) \ = \ \inf_{\beta[\cdot]\in\Bc_{\textup{\tiny str}}} \sup_{\alpha\in\Ac} V\big(s,X_s^{t,\xi,\alpha,\beta[\alpha]},\tilde X_s^{t,\xi,\eta,\alpha,\beta[\alpha]}\big),
\end{equation}
for all $t,s\in[0,T]$, with $t\leq s$, and for every $(\xi,\eta)\in L^q(\Omega,\Fc_t\vee\Gc,\P;\R^n)\times L^1(\Omega,\Fc_t\vee\Gc,\P;\R)$. Hence, it remains to prove \eqref{DPPtildeV}. The following proof of \eqref{DPPtildeV} is inspired by the proof of the dynamic programming principle for deterministic differential games, see Theorem 3.1 in \cite{EvansSouganidis}.

Fix $t,s\in[0,T]$, with $t\leq s$, and $(\xi,\eta)\in L^q(\Omega,\Fc_t\vee\Gc,\P;\R^n)\times L^1(\Omega,\Fc_t\vee\Gc,\P;\R)$. Set
\[
\Lambda(t,\xi,\eta) \ := \ \inf_{\beta[\cdot]\in\Bc_{\textup{\tiny str}}} \sup_{\alpha\in\Ac} V\big(s,X_s^{t,\xi,\alpha,\beta[\alpha]},\tilde X_s^{t,\xi,\eta,\alpha,\beta[\alpha]}\big).
\]
We split the proof of \eqref{DPPtildeV} into two steps.

\vspace{2mm}

\noindent\emph{Proof of $V(t,\xi,\eta)\leq\Lambda(t,\xi,\eta)$.} Fix $\eps>0$. Then, there exists ${\bar\beta}^\eps[\cdot]\in\Bc_{\textup{\tiny str}}$ such that
\begin{align}\label{Ineq1}
\Lambda(t,\xi,\eta) \ &\geq \ \sup_{\alpha\in\Ac} V\big(s,X_s^{t,\xi,\alpha,{\bar\beta}^\eps[\alpha]},\tilde X_s^{t,\xi,\eta,\alpha,{\bar\beta}^\eps[\alpha]}\big) - \eps \notag \\
&\geq \ V\big(s,X_s^{t,\xi,\alpha,{\bar\beta}^\eps[\alpha]},\tilde X_s^{t,\xi,\eta,\alpha,{\bar\beta}^\eps[\alpha]}\big) - \eps, \qquad \text{for every }\alpha\in\Ac.
\end{align}
Now, notice that for every fixed $\alpha\in\Ac$ there exists $\bar{\beta'}^{,\eps,\alpha}[\cdot]\in\Bc_{\textup{\tiny str}}$ such that
\begin{align}\label{Ineq2}
&V\big(s,X_s^{t,\xi,\alpha,{\bar\beta}^\eps[\alpha]},\tilde X_s^{t,\xi,\eta,\alpha,{\bar\beta}^\eps[\alpha]}\big) \notag \\
&= \ \inf_{\beta'[\cdot]\in\Bc_{\textup{\tiny str}}} \sup_{\alpha'\in\Ac} G\Big(X_T^{s,X_s^{t,\xi,\alpha,\bar\beta^\eps[\alpha]},\alpha',\beta'[\alpha']},\tilde X_T^{s,X_s^{t,\xi,\alpha,\bar\beta^\eps[\alpha]},\tilde X_s^{t,\xi,\eta,\alpha,\bar\beta^\eps[\alpha]},\alpha',\beta'[\alpha']}\Big) \notag \\
&\geq \ \sup_{\alpha'\in\Ac} G\Big(X_T^{s,X_s^{t,\xi,\alpha,\bar\beta^\eps[\alpha]},\alpha',\bar{\beta'}^{,\eps,\alpha}[\alpha']},\tilde X_T^{s,X_s^{t,\xi,\alpha,\bar\beta^\eps[\alpha]},\tilde X_s^{t,\xi,\eta,\alpha,\bar\beta^\eps[\alpha]},\alpha',\bar{\beta'}^{,\eps,\alpha}[\alpha']}\Big) - \eps \notag \\
&\geq \ G\Big(X_T^{s,X_s^{t,\xi,\alpha,\bar\beta^\eps[\alpha]},\alpha,\bar{\beta'}^{,\eps,\alpha}[\alpha]},\tilde X_T^{s,X_s^{t,\xi,\alpha,\bar\beta^\eps[\alpha]},\tilde X_s^{t,\xi,\eta,\alpha,\bar\beta^\eps[\alpha]},\alpha,\bar{\beta'}^{,\eps,\alpha}[\alpha]}\Big) - \eps.
\end{align}
Define $\beta^\eps[\cdot]\in\Bc_{\textup{\tiny str}}$ as follows: for every fixed $\alpha\in\Ac$, we set
\[
\beta^\eps[\alpha]_r \ := \ \bar\beta^\eps[\alpha]_r\,1_{[0,s]}(r) + \bar{\beta'}^{,\eps,\alpha}[\alpha]_r\,1_{(s,T]}(r), \qquad \text{for all }r\in[0,T].
\]
Then, we can rewrite \eqref{Ineq2} in terms of $\beta^\eps$ as
\begin{align*}
&V\big(s,X_s^{t,\xi,\alpha,{\bar\beta}^\eps[\alpha]},\tilde X_s^{t,\xi,\eta,\alpha,{\bar\beta}^\eps[\alpha]}\big) \\
&= V\big(s,X_s^{t,\xi,\alpha,\beta^\eps[\alpha]},\tilde X_s^{t,\xi,\eta,\alpha,\beta^\eps[\alpha]}\big) \geq G\Big(X_T^{s,X_s^{t,\xi,\alpha,\beta^\eps[\alpha]},\alpha,\beta^\eps[\alpha]},\tilde X_T^{s,X_s^{t,\xi,\alpha,\beta^\eps[\alpha]},\tilde X_s^{t,\xi,\eta,\alpha,\beta^\eps[\alpha]},\alpha,\beta^\eps[\alpha]}\Big) - \eps.
\end{align*}
By the flow properties \eqref{flow}-\eqref{flow_tildeX}, we obtain
\[
V\big(s,X_s^{t,\xi,\alpha,\beta^\eps[\alpha]},\tilde X_s^{t,\xi,\eta,\alpha,\beta^\eps[\alpha]}\big) \ \geq \ G\big(X_T^{t,\xi,\alpha,\beta^\eps[\alpha]},\tilde X_s^{t,\xi,\eta,\alpha,\beta^\eps[\alpha]}\big) - \eps.
\]
Plugging the above inequality into \eqref{Ineq1}, we find
\[
\Lambda(t,\xi,\eta) \ \geq \ G\big(X_T^{t,\xi,\alpha,\beta^\eps[\alpha]},\tilde X_s^{t,\xi,\eta,\alpha,\beta^\eps[\alpha]}\big) - 2\,\eps, \qquad \text{for every }\alpha\in\Ac.
\]
The claim follows taking the supremum over $\Ac$ and then the infimum over $\Bc_{\textup{\tiny str}}$.

\vspace{2mm}

\noindent\emph{Proof of $V(t,\xi,\eta)\geq\Lambda(t,\xi,\eta)$.} Fix $\eps>0$. Then, there exists ${\bar\beta}^\eps[\cdot]\in\Bc_{\textup{\tiny str}}$ such that
\begin{equation}\label{Ineq3}
\sup_{\alpha\in\Ac} G\big(X_T^{t,\xi,\alpha,{\bar\beta}^\eps[\alpha]},\tilde X_T^{t,\xi,\eta,\alpha,{\bar\beta}^\eps[\alpha]}\big) \ \leq \ V(t,\xi,\eta) + \eps.
\end{equation}
We also have
\[
\Lambda(t,\xi,\eta) \ \leq \ \sup_{\alpha\in\Ac} V\big(s,X_s^{t,\xi,\alpha,{\bar\beta}^\eps[\alpha]},\tilde X_s^{t,\xi,\eta,\alpha,{\bar\beta}^\eps[\alpha]}\big).
\]
So, in particular, there exists $\alpha^\eps\in\Ac$ such that
\begin{equation}\label{Ineq4}
\Lambda(t,\xi,\eta) \ \leq \ V\big(s,X_s^{t,\xi,\alpha^\eps,{\bar\beta}^\eps[\alpha^\eps]},\tilde X_s^{t,\xi,\eta,\alpha^\eps,{\bar\beta}^\eps[\alpha^\eps]}\big) + \eps.
\end{equation}
Now, for every $\alpha\in\Ac$ define $\tilde\alpha^\eps\in\Ac$ by
\begin{equation}\label{tilde_alpha}
\tilde\alpha^\eps_r \ := \ \alpha_r^\eps\,1_{[0,s]}(r) + \alpha_r\,1_{(s,T]}(r), \qquad \text{for all }r\in[0,T].
\end{equation}
Then, define $\beta^\eps\in\Bc_{\textup{\tiny str}}$ by
$\beta^\eps[\alpha] \ := \ \bar\beta^\eps[\tilde\alpha^\eps]$,  for every $\alpha$ $\in$ $\Ac$. 
Hence
\begin{align*}
&V\big(s,X_s^{t,\xi,\alpha^\eps,{\bar\beta}^\eps[\alpha^\eps]},\tilde X_s^{t,\xi,\eta,\alpha^\eps,{\bar\beta}^\eps[\alpha^\eps]}\big) \\
&= \ \inf_{\beta[\cdot]\in\Bc_{\textup{\tiny str}}} \sup_{\alpha\in\Ac} G\Big(X_T^{s,X_s^{t,\xi,\alpha^\eps,\bar\beta^\eps[\alpha^\eps]},\alpha,\beta[\alpha]},\tilde X_T^{s,X_s^{t,\xi,\alpha^\eps,\bar\beta^\eps[\alpha^\eps]},\tilde X_s^{t,\xi,\eta,\alpha^\eps,\bar\beta^\eps[\alpha^\eps]},\alpha,\beta[\alpha]}\Big) \\
&\leq \ \sup_{\alpha\in\Ac} G\Big(X_T^{s,X_s^{t,\xi,\alpha^\eps,\bar\beta^\eps[\alpha^\eps]},\alpha,\beta^\eps[\alpha]},\tilde X_T^{s,X_s^{t,\xi,\alpha^\eps,\bar\beta^\eps[\alpha^\eps]},\tilde X_s^{t,\xi,\eta,\alpha^\eps,\bar\beta^\eps[\alpha^\eps]},\alpha,\beta^\eps[\alpha]}\Big) \\
&= \ \sup_{\alpha\in\Ac} G\Big(X_T^{s,X_s^{t,\xi,\tilde\alpha^\eps,\bar\beta^\eps[\tilde\alpha^\eps]},\tilde\alpha^\eps,\bar\beta^\eps[\tilde\alpha^\eps]},\tilde X_T^{s,X_s^{t,\xi,\tilde\alpha^\eps,\bar\beta^\eps[\tilde\alpha^\eps]},\tilde X_s^{t,\xi,\eta,\tilde\alpha^\eps,\bar\beta^\eps[\tilde\alpha^\eps]},\tilde\alpha^\eps,\bar\beta^\eps[\tilde\alpha^\eps]}\Big),
\end{align*}
where the last equality follows from the definitions of $\tilde\alpha^\eps$ and $\beta^\eps$. By the flow properties \eqref{flow}-\eqref{flow_tildeX}, we obtain
\[
V\big(s,X_s^{t,\xi,\alpha^\eps,{\bar\beta}^\eps[\alpha^\eps]},\tilde X_s^{t,\xi,\eta,\alpha^\eps,{\bar\beta}^\eps[\alpha^\eps]}\big) \ \leq \ \sup_{\alpha\in\Ac} G\big(X_T^{t,\xi,\tilde\alpha^\eps,\bar\beta^\eps[\tilde\alpha^\eps]},\tilde X_T^{t,\xi,\eta,\tilde\alpha^\eps,\bar\beta^\eps[\tilde\alpha^\eps]}\big).
\]
Consequently, there exists $\alpha^{2,\eps}\in\Ac$, and the corresponding $\tilde\alpha^{2,\eps}$ defined as in \eqref{tilde_alpha}, such that
\[
V\big(s,X_s^{t,\xi,\alpha^\eps,{\bar\beta}^\eps[\alpha^\eps]},\tilde X_s^{t,\xi,\eta,\alpha^\eps,{\bar\beta}^\eps[\alpha^\eps]}\big) \ \leq \ G\big(X_T^{t,\xi,\tilde\alpha^{2,\eps},\bar\beta^\eps[\tilde\alpha^{2,\eps}]},\tilde X_T^{t,\xi,\eta,\tilde\alpha^{2,\eps},\bar\beta^\eps[\tilde\alpha^{2,\eps}]}\big) + \eps.
\]
Finally, using inequalities \eqref{Ineq3} and \eqref{Ineq4}, we obtain
\begin{align*}
\Lambda(t,\xi,\eta) \ &\leq \ V\big(s,X_s^{t,\xi,\alpha^\eps,{\bar\beta}^\eps[\alpha^\eps]},\tilde X_s^{t,\xi,\eta,\alpha^\eps,{\bar\beta}^\eps[\alpha^\eps]}\big) + \eps \\
&\leq \ G\big(X_T^{t,\xi,\tilde\alpha^{2,\eps},\bar\beta^\eps[\tilde\alpha^{2,\eps}]},\tilde X_T^{t,\xi,\eta,\tilde\alpha^{2,\eps},\bar\beta^\eps[\tilde\alpha^{2,\eps}]}\big) + 2\,\eps \\
&\leq \ \sup_{\alpha\in\Ac} G\big(X_T^{t,\xi,\alpha,\bar\beta^\eps[\alpha]},\tilde X_T^{t,\xi,\eta,\alpha,\bar\beta^\eps[\alpha]}\big) + 2\,\eps \ \leq \ V(t,\xi,\eta) + 3\,\eps,
\end{align*}
which concludes the proof.
\ep

\vspace{3mm}

We immediately deduce the DPP for the inverse-lifted  lower and upper value functions. 

\begin{Corollary}\label{invlift:DPP}
Under Assumption {\bf (A1)}, we have
\begin{align}\label{DPPvinv}
\mathscr V(t,\mu) \ = \ \inf_{\beta[\cdot]\in\Bc_{\textup{\tiny str}}} \sup_{\alpha\in\Ac} \bigg\{ \E\Big[\int_t^s f\big(X_r^{t,\xi,\alpha,\beta[\alpha]},\P_{\text{\tiny$X_r^{t,\xi,\alpha,\beta[\alpha]}$}},\alpha_r,\beta[\alpha]_r,\P_{\text{\tiny$(\alpha_r,\beta[\alpha]_r)$}}\big)\,dr&\Big] \\
\;\;\; +  \;\;  \mathscr V\big(s,\P_{\text{\tiny$X_s^{t,\xi,\alpha,\beta[\alpha]}$}}\big)&\bigg\} \notag
\end{align}
and
\begin{align*}
\mathscr U(t,\mu)\ = \ \sup_{\alpha[\cdot]\in\Ac_{\textup{\tiny str}}} \inf_{\beta\in\Bc} \bigg\{ \E\Big[\int_t^s f\big(X_r^{t,\xi,\alpha[\beta],\beta},\P_{\text{\tiny$X_r^{t,\xi,\alpha[\beta],\beta}$}},\alpha[\beta]_r,\beta_r,\P_{\text{\tiny$(\alpha[\beta]_r,\beta_r)$}}\big)\,dr&\Big] \\
+ \, \mathscr U\big(s,\P_{\text{\tiny$X_s^{t,\xi,\alpha,\beta[\alpha]}$}}\big)&\bigg\}, \notag
\end{align*}
for all $t,s\in[0,T]$, with $t\leq s$, $\mu$ $\in$ $\mathscr P_{\text{\tiny$q$}}(\R^n)$, and any $\xi\in L^q(\Omega,\Gc,\P;\R^n)$ such that $\P_{\text{\tiny$\xi$}}$ $=$ $\mu$. 
\end{Corollary}

\vspace{5mm}

\noindent {\bf The case without mean-field interaction}

\vspace{2mm}

\noindent Let us consider the particular case of standard stochastic optimal control problem (for the case of a standard two-player zero-sum stochastic differential game see Remark \ref{R:ClassicalGame} below), where $\Gc$ is the trivial $\sigma$-algebra and all coefficients depend only on the state and control (of the first player), but not on their probability laws (as well as on the control of the second player):

\vspace{2mm}

\noindent {\bf (A3)}\;\;$\Gc=\{\emptyset,\Omega\}$ and $\gamma$ $=$ $\gamma(x,a)$, $\sigma$ $=$ $\sigma(x,a)$, $f$ $=$ $f(x,a)$, $g$ $=$ $g(x)$, for every $(x,a)$ $\in$ $\R^n\times A$.

\vspace{2mm}

The assumption that $\Gc$ is trivial implies that the lower and upper value functions $v=v(t,\xi)$ and $u=u(t,\xi)$ are defined only for all $t\in[0,T]$, $\xi\in L^q(\Omega,\Fc_t,\P;\R^n)$ (rather than $\xi\in L^q(\Omega,\Fc_t\vee\Gc,\P;\R^n)$). This also reflects on $\mathscr V$ and $\mathscr U$, which now are defined on a possibly smaller set, given by all pairs $(t,\mu)\in[0,T]\times\mathscr P_{\text{\tiny$q$}}(\R^n)$ for which there exists $\xi\in L^q(\Omega,\Fc_t,\P;\R^n)$ such that $\P_{\text{\tiny$\xi$}}=\mu$. Another consequence of the assumption that $\Gc$ is trivial is that $\Ac$ (resp. $\Bc$) coincides with the family of all $\F$-progressively (rather than $(\Fc_s\vee\Gc)_s$-progressively) measurable processes taking values in $A$ (resp. $B$).

Under Assumption {\bf (A3)}, denote by $X^{t,\xi,\alpha}$ the solution to \reff{State} when there is only the control process $\alpha$. By an abuse of notation, we still denote by $J=J(t,\xi,\alpha)$ the payoff functional, which now depends only on $t$, $\xi$, $\alpha$:
\[
J(t,\xi,\alpha) \ = \ \E\bigg[\int_t^T f\big(X_s^{t,\xi,\alpha},\alpha_s\big)\,ds + g\big(X_T^{t,\xi,\alpha}\big)\bigg].
\]
Notice that under {\bf (A3)}, the lower and upper value functions coincide with each other, and are simply given by:
\[
v(t,\xi) \ = \ u(t,\xi) \ = \ \sup_{\alpha\in\Ac} J(t,\xi,\alpha), \qquad \text{for all }t\in[0,T],\,\xi\in L^q(\Omega,\Fc_t,\P;\R^n).
\]
For any $t$ $\in$ $[0,T]$, $x$ $\in$ $\R^n$, $\alpha$ $\in$ $\Ac$, denote by $X^{t,x,\alpha}$ the solution to \reff{State} when the initial condition at time $t$ is given by a constant $\xi$ $=$ $x$ in $\R^n$. Similarly, let $J^B(t,x,\alpha)$ be the payoff functional when $\xi$ $=$ $x$, namely (the capital letter $B$ at the top of $J$ refers to ``Bellman'') 
\[
J^B(t,x,\alpha) \ = \ \E\bigg[\int_t^T f\big(X_s^{t,x,\alpha},\alpha_s\big)\,ds + g\big(X_T^{t,x,\alpha}\big)\bigg].
\]
So, in particular, $J^B(t,x,\alpha)$ coincides with $J(t,\xi,\alpha)$ whenever $\xi$ $=$ $x$ in $\R^n$. Finally, the value function of the standard stochastic optimal control problem is given by:
\[
v^B(t,x) \ = \ \sup_{\alpha\in\Ac} J^B(t,x,\alpha),
\]
for all $t$ $\in$ $[0,T]$, $x$ $\in$ $\R^n$.
 
The following result makes the connection between the standard value function $v^B$ and our value functions $v=u$, or the value functions $\mathscr V=\mathscr U$ on the Wasserstein space, and show that one can retrieve the standard DPP in the non-McKean-Vlasov case from Theorem \ref{T:DPP}. 

\begin{Proposition}\label{P:Classical}
Under Assumptions {\bf (A1)} and {\bf (A3)}, we have, for all $t$ $\in$ $[0,T]$, $\mu$ $\in$ $\mathscr P_{\text{\tiny$q$}}(\R^n)$,
\begin{equation}\label{Identities}
\mathscr V(t,\mu) \; = \; v(t,\xi) \; = \; u(t,\xi) \; = \; \mathscr U(t,\mu) \; = \; \int_{\R^n} v^B(t,x) \mu(dx) \; = \; \E[v^B(t,\xi)],
\end{equation}
for any $\xi\in L^q(\Omega,\Fc_t,\P;\R^n)$, with $\P_{\text{\tiny$\xi$}}=\mu$. Therefore, we have the DPP for $v^B$:
\beqs
v^B(t,x) & = & \sup_{\alpha\in\Ac}  \E\bigg[\int_t^s f\big(X_r^{t,x,\alpha},\alpha_r\big)\,dr   +   
v^B\big(s,X_s^{t,x,\alpha}\big) \bigg],
\enqs
for all $t,s\in[0,T]$, with $t\leq s$, and for every $x\in\R^n$.
\end{Proposition}
{\bf Proof.}
If the sequence of equalities \eqref{Identities} holds true, the DPP for $v^B$ is a direct consequence of Theorem \ref{T:DPP}. Then, it remains to prove \eqref{Identities}.

Fix $t\in[0,T]$, $\mu\in\mathscr P_{\text{\tiny$q$}}(\R^n)$, and $\xi\in L^q(\Omega,\Fc_t,\P;\R^n)$ with $\P_{\text{\tiny$\xi$}}=\mu$. Equalities $\mathscr V=v=u=\mathscr U$ follow directly from Assumption {\bf (A3)}. Therefore, it only remains to prove the following equality:
\[
v(t,\xi) \ = \ \E[v^B(t,\xi)].
\]

\vspace{2mm}

\noindent\emph{\textbf{Proof of the inequality $v(t,\xi)$ $\leq$ $\E[v^B(t,\xi)]$.}} We adopt the same notations as in \textbf{Steps I} and \textbf{II} of the proof of Proposition \ref{P:Lift}. So, in particular, we consider a uniformly distributed random variable $\Gamma_t\colon\Omega\rightarrow[0,1]$ such that $\sigma(\Gamma_t)=\Fc_t^o$. Moreover
\begin{align}
\xi \ &= \ \chi(\Gamma_t), \qquad\qquad \P\text{-a.s.}, \notag \\
\alpha_s \ &= \ a_s(\cdot,\Gamma_t(\cdot)), \qquad \forall\,s\in[0,T],\,\P\text{-a.s.}, \label{alpha=a}
\end{align}
for some measurable functions $\chi\colon([0,1],\Bc([0,1]))\rightarrow(\R^n,\Bc(\R^n))$ and $a\colon(\Omega\times[0,T]\times[0,1],$ $Prog(\F^t)\otimes\Bc([0,1]))\rightarrow(A,\Bc(A))$. For every fixed $y\in[0,1]$, denote $\alpha^y:=a_\cdot(\cdot,y)$. Notice that $\alpha^y\in\Ac$ and, in particular, it is $\F^t$-progressively measurable.

Similarly, the controlled state process $(X_s^{t,\xi,\alpha})_{s\in[t,T]}$ is $\F$-progressively measurable, so, in particular, it is $Prog(\F^t)\vee(\Fc_t\otimes\{\emptyset,[0,T]\})$-measurable (where $\{\emptyset,[0,T]\}$ denotes the trivial $\sigma$-algebra on $[0,T]$). Then, $X^{t,\xi,\alpha}$ has the form
\[ 
X_s^{t,\xi,\alpha} \ = \ \boldsymbol x_s^{t,\xi,\alpha}(\cdot,\Gamma_t(\cdot)), \qquad \forall\,s\in[0,T],\,\P\text{-a.s.},
\]
for some $Prog(\F^t)\otimes\Bc([0,1])$-measurable function $\boldsymbol x^{t,\xi,\alpha}=\boldsymbol x_s^{t,\xi,\alpha}(\omega,y)\colon\Omega\times[0,T]\times[0,1]\rightarrow\R^n$. By pathwise uniqueness to equation \eqref{State}, we deduce that there exists a Lebesgue-null set $N\in\Bc([0,1])$ such that
\[
\boldsymbol x_s^{t,\xi,\alpha}(\cdot,y) \ = \ X_s^{t,\chi(y),\alpha^y}, \qquad \forall\,s\in[t,T],\,\P\text{-a.s.},
\]
for every $y\in[0,1]\backslash N$. Hence, using the fact that $\sigma(\Gamma_t)=\Fc_t^o$ is independent of $\Fc_T^t$, by Fubini's theorem we deduce that the payoff functional can be written as follows:
\begin{align}
J(t,\xi,\alpha) \ &= \ \int_0^1\E\bigg[\int_t^T f\big(X_s^{t,\chi(y),\alpha^y},\alpha_s^y\big)\,ds + g\big(X_T^{t,\chi(y),\alpha^y}\big)\bigg]dy \label{J=JBy} \\
&= \ \int_0^1 J^B\big(t,\chi(y),\alpha^y\big)\,dy \ \leq \ \int_0^1 v^B(t,\chi(y))\,dy. \notag
\end{align}
Taking the supremum over $\alpha$ in $\Ac$, we conclude that
\[
v(t,\xi) \ \leq \ \int_0^1 v^B(t,\chi(y))\,dy \ = \ \E[v^B(t,\xi)].
\]

\vspace{2mm}

\noindent\emph{\textbf{Proof of the inequality $v(t,\xi)$ $\geq$ $\E[v^B(t,\xi)]$.}} Recall that
\[
\E[v^B(t,\xi)] \ = \ \int_{\R^n} v^B(t,x)\,\mu(dx)
\]
and, for every $x\in\R^n$,
\[
v^B(t,x) \ = \ \sup_{\alpha\in\Ac} J^B(t,x,\alpha) \ = \ \sup_{\alpha\in\Ac} \E\bigg[\int_t^T f\big(X_s^{t,x,\alpha},\alpha_s\big)\,ds + g\big(X_T^{t,x,\alpha}\big)\bigg].
\]
Notice also that the following equality holds:
\begin{equation}\label{Ac^t}
\sup_{\alpha\in\Ac} J^B(t,x,\alpha) \ = \ \sup_{\alpha\in\Ac^t} J^B(t,x,\alpha),
\end{equation}
with $\Ac^t\subset\Ac$ denoting the set of all $\F^t$-progressively measurable processes taking values in $A$. In order to see that \eqref{Ac^t} holds, we begin noting that $\sup_{\alpha\in\Ac} J^B(t,x,\alpha)\geq\sup_{\alpha\in\Ac^t} J^B(t,x,\alpha)$ since $\Ac^t\subset\Ac$, so it remains to prove the other inequality. Given $\alpha\in\Ac$, recalling that $\alpha$ can be written as in \eqref{alpha=a}, we find (proceeding as in \eqref{J=JBy})
\begin{align*}
J^B(t,x,\alpha) \ &= \ \int_0^1\E\bigg[\int_t^T f\big(X_s^{t,x,\alpha^y},\alpha_s^y\big)\,ds + g\big(X_T^{t,x,\alpha^y}\big)\bigg]dy \\
&= \ \int_0^1 J^B(t,x,\alpha^y)\,dy \ \leq \ \int_0^1\sup_{\tilde\alpha\in\Ac^t}J^B(t,x,\tilde\alpha)\,dy \ = \ \sup_{\tilde\alpha\in\Ac^t}J^B(t,x,\tilde\alpha),
\end{align*}
where the inequality follows from the fact that $\alpha^y\in\Ac^t$, for every $y\in[0,1]$. From the arbitrariness of $\alpha\in\Ac$, we deduce the other inequality, from which \eqref{Ac^t} follows.

Given $\eps>0$, for every $x\in\R^n$ let $\alpha^{x,\eps}\in\Ac^t$ be an $\eps$-optimal control, namely
\begin{equation}\label{alpha^x,eps}
v^B(t,x) \ \leq \ J^B(t,x,\alpha^{x,\eps}) + \eps.
\end{equation}
Suppose for a moment that the composition $\alpha^\eps:=\alpha^{\xi,\eps}$ belongs to $\Ac$. Then, from  \eqref{alpha^x,eps} we get
\begin{equation}\label{alpha^x,eps2}
v^B(t,\xi) \ \leq \ J^B(t,\xi,\alpha^\eps) + \eps.
\end{equation}
Notice that $J^B(t,\xi,\alpha^\eps)$ is a random variable, indeed it is a function of $\xi$. We also observe, similarly  as in \reff{J=JBy}, that the expectation of $J^B(t,\xi,\alpha^\eps)$ coincides with $J(t,\xi,\alpha^\eps)$. Then, taking the expectation in \eqref{alpha^x,eps2}, we obtain
\[
\E[v^B(t,\xi)] \ \leq \ \E[J^B(t,\xi,\alpha^\eps)] + \eps \ = \ J(t,\xi,\alpha^\eps) + \eps \ \leq \ v(t,\xi) + \eps.
\]
From the arbitrariness of $\eps$, the claim follows. It remains to prove that $\alpha^{\xi,\eps}$ belongs to $\Ac$. More precisely, it is enough to prove that for every $x\in\R^n$ we are able to find an $\eps$-optimal control $\alpha^{x,\eps}\in\Ac^t$ such that the composition $\alpha^{\xi,\eps}$ belongs to $\Ac$. This last statement follows easily when the random variable $\xi$ is discrete. In the general case, we need to apply a measurable selection theorem. To this end, we define the following metric on $\Ac^t$ (see Definition 3.2.3 in \cite{Krylov}):
\[
\rho_{\textup{Kr}}(\alpha,\alpha') \ := \ \E\bigg[\int_0^T \rho_A(\alpha_s,\alpha_s')\,ds\bigg],
\]
for any $\alpha,\alpha'\in\Ac^t$, where we recall that $\rho_A$ is a bounded metric on the Polish space $A$. We observe that the metric space $(\Ac^t,\rho_{\textup{Kr}})$ is complete. Let us now prove that $(\Ac^t,\rho_{\textup{Kr}})$ is also separable, so $(\Ac^t,\rho_{\textup{Kr}})$ is a Polish space. Firstly, notice that if $A$ is equal to some Euclidean space $\R^m$, then $\Ac^t$ coincides with the closed subset of the space $L^1(\Omega\times[0,T],\Fc_T\otimes\Bc([0,T]),d\P\otimes ds;(A,\rho_A))$ of all $\F^t$-progressively measurable processes; so, in particular, (see page 92 in \cite{Doob} and the beginning of Section 2.5 in \cite{SonerTouzi}) the space $(\Ac^t,\rho_{\textup{Kr}})$ is separable since the $\sigma$-algebra $\Fc_T^t\otimes\Bc([0,T])$ is countably generated up to null sets. When $A$ is a generic Polish space, the same result holds true. As a matter fact, proceeding as in \cite{Doob}, page 92, we see that the separability of $(\Ac^t,\rho_{\textup{Kr}})$ follows from the following facts: since, up to null sets, the $\sigma$-algebra $\Fc_T^t\otimes\Bc([0,T])$ is countably generated, the subfamily of $\Ac^t$ of all processes which are equal to an indicator function of a measurable set in $\Fc_T^t\otimes\Bc([0,T])$ is separable; it follows that the subfamily of all step processes taking values in some fixed countable subset of $A$ is separable; given that $A$ is separable, the subfamily of all step processes is separable; since this latter subfamily is dense in $\Ac^t$, it follows that $\Ac^t$ is also separable.

We can now apply the Jankov-von Neumann measurable selection theorem, and in particular Proposition 7.50 in \cite{BertsekasShreve}. More precisely, in order to apply Proposition 7.50 in \cite{BertsekasShreve}, we begin noting that $X$, $Y$, $D$, $f$, $f^*$ in \cite{BertsekasShreve} are given respectively by $\R^n$, $\Ac^t$, $\R^n\times\Ac^t$, $-J^B(t,\cdot,\cdot)$ (the minus sign is due to the presence of the \emph{inf} in \cite{BertsekasShreve}), $v^B(t,\cdot)$. We firstly notice that $-J^B(t,\cdot,\cdot)\colon\R^n\times\Ac^t\rightarrow\R$ is a Borel measurable function, so, in particular, it is lower semianalytic. Then, by Proposition 7.50 in \cite{BertsekasShreve} it follows that: for any $\eps>0$, there exists an analytically measurable function $\boldsymbol\alpha^\eps\colon\R^n\rightarrow\Ac^t$ such that
\[
v^B(t,x) \ \leq \ J^B(t,x,\boldsymbol\alpha^\eps(x)) + \eps, \qquad \text{for every }x\in\R^n.
\]
Since every analytic set in $\R^n$ belongs to the Lebesgue $\sigma$-algebra $\mathscr L(\R^n)$, we see that $\boldsymbol\alpha^\eps$ is a measurable function from $(\R^n,\mathscr L(\R^n))$ into $(\Ac^t,\Bc(\Ac^t))$. Now, it is easy to see that there exists a measurable function $\tilde{\boldsymbol\alpha}^\eps\colon(\R^n,\Bc(\R^n))\rightarrow(\Ac^t,\Bc(\Ac^t))$ which is equal to $\boldsymbol\alpha^\eps$ a.e. (with respect to the Lebesgue measure on $\R^n$). As a matter of fact, if $\boldsymbol\alpha^\eps$ is a sum of indicator functions on a Lebesgue measurable partition of $\R^n$, the result follows from the fact that if $\mathscr B\in\mathscr L(\R^n)$ then there exists $\tilde{\mathscr B}\in\Bc(\R^n)$ such that the Lebesgue measure of $\mathscr B\Delta\tilde{\mathscr B}$ is zero; for a general Lebesgue measurable function $\boldsymbol\alpha^\eps$, the result follows by an approximation argument. We thus have
\[
v^B(t,x) \ \leq \ J^B(t,x,\tilde{\boldsymbol\alpha}^\eps(x)) + \eps, \qquad \text{for a.e. }x\in\R^n.
\]
In order to conclude the proof, we notice that it remains to prove that the composition $\alpha^\eps:=\tilde{\boldsymbol\alpha}^\eps(\xi)$ belongs to $\Ac$. To this end, suppose firstly that $\tilde{\boldsymbol\alpha}^\eps$ is a sum of indicator functions on a Borel measurable partition of $\R^n$, namely
\begin{equation}\label{tildealpha^eps}
\tilde{\boldsymbol\alpha}^\eps(x) \ = \ \sum_i \alpha_i^\eps\,1_{\{x\in\mathscr B_i\}},
\end{equation}
where $\{\alpha_i^\eps\}_i\subset\Ac^t$ and $\{\mathscr B_i\}_i\subset\Bc(\R^n)$ is a partition of $\R^n$. If $\tilde{\boldsymbol\alpha}^\eps$ has the form in \eqref{tildealpha^eps}, then it is clear that the composition $\tilde{\boldsymbol\alpha}^\eps(\xi)$ belongs to $\Ac$, that is $\tilde{\boldsymbol\alpha}^\eps(\xi)$ is an $\F$-progressively measurable process (as a matter of fact, for every $i$, both $\alpha_i^\eps$ and the indicator function $1_{\{\xi\in\mathscr B_i\}}$ are $\F$-progressively measurable processes). For a general Borel measurable function, the result follows by an approximation argument.
\ep

\begin{Remark}\label{R:ClassicalGame}
{\rm
The extension of Proposition \ref{P:Classical} to the case of two-player zero-sum stochastic differential games presents some difficulties, as we now explain. Consider a standard two-player zero-sum stochastic differential game, firstly studied in the seminal paper \cite{FlemingSouganidis}, where  all coefficients depend only on the state and controls, but not on their probability laws:

\vspace{1mm}

\noindent {\bf (A3)}${}_{\text{game}}$\;\;$\Gc=\{\emptyset,\Omega\}$ and $\gamma$ $=$ $\gamma(x,a,b)$, $\sigma$ $=$ $\sigma(x,a,b)$, $f$ $=$ $f(x,a,b)$, $g$ $=$ $g(x)$, for every $(x,a,b)$ $\in$ $\R^d\times A\times B$

\vspace{1mm}

The lower and upper value functions of this standard stochastic differential game, as consi\-dered in \cite{FlemingSouganidis}, are defined as follows:
\beqs
v^{FS}(t,x) & = &  \inf_{\beta[\cdot]\in\Bc_{\textup{\tiny str}}} \sup_{\alpha\in\Ac} J^{FS}(t,x,\alpha,\beta[\alpha]), \\ 
u^{FS}(t,x) & = & \sup_{\alpha[\cdot]\in\Ac_{\textup{\tiny str}}} \inf_{\beta\in\Bc} J^{FS}(t,x,\alpha[\beta],\beta),
\enqs
for all $t\in[0,T]$, $x\in\R^n$, where
\[
J^{FS}(t,x,\alpha,\beta) \ := \ \E\bigg[\int_t^T f\big(s,X_s^{t,x,\alpha,\beta},\alpha_s,\beta_s\big)\,ds + g\big(X_T^{t,x,\alpha,\beta}\big)\bigg].
\]
In order to adapt the proof of Proposition \ref{P:Classical} to this context, we need the following generalizations of equality \eqref{Ac^t}:
\begin{align}
\inf_{\beta[\cdot]\in\Bc_{\textup{\tiny str}}} \sup_{\alpha\in\Ac} 
J^{FS}(t,x,\alpha,\beta[\alpha]) \ &= \ \inf_{\beta[\cdot]\in\Bc_{\textup{\tiny str}}^t} \sup_{\alpha\in\Ac^t} 
J^{FS}(t,x,\alpha,\beta[\alpha]), \label{Ac^t1} \\
 \sup_{\alpha[\cdot]\in\Ac_{\textup{\tiny str}}} \inf_{\beta\in\Bc} J^{FS}(t,x,\alpha[\beta],\beta) \ &= \  \sup_{\alpha[\cdot]\in\Ac_{\textup{\tiny str}}^t} \inf_{\beta\in\Bc^t} J^{FS}(t,x,\alpha[\beta],\beta), \label{Ac^t2}
\end{align}
where $\Ac^t$ (resp. $\Bc^t$) denotes the set of $\F^t$-progressively measurable processes taking values in $A$ (resp. $B$), while $\Ac_{\textup{\tiny str}}^t$ (resp. $\Bc_{\textup{\tiny str}}^t$) is defined as $\Ac_{\textup{\tiny str}}$ (resp. $\Bc_{\textup{\tiny str}}$) replacing $\Ac$ and $\Bc$ respectively by $\Ac^t$ and $\Bc^t$. The validity of \eqref{Ac^t1} and \eqref{Ac^t2} is however, at least to our knowledge, still not known in the literature, 
although we conjecture that relations: $v(t,\xi)$ $=$ $\E[v^{FS}(t,\xi)]$, $u(t,\xi)$ $=$ $\E[u^{FS}(t,\xi)]$ hold true, see also Remark \ref{remHJBI}. 
\epR}
\end{Remark}

\section{Bellman-Isaacs dynamic programming equations} \label{secvisco} 

This section is devoted to the derivation of the Bellman-Isaacs equation from the DPP for the lower and upper value functions, and the viscosity PDE (partial differential equation) characterization. 
We shall provide a PDE formulation for the value functions on the Hilbert space $L^2(\Omega,\Gc,\P;\R^n)$  (hence for $q$ $=$ $2$) or alternatively via the inverse-lifted (recall Definition \ref{definverselift}) identification on the Wasserstein space of probability measures $\mathscr P_{\text{\tiny$2$}}(\R^n)$. 

We shall rely on the notion of lifted derivative with respect to a probability measure and It\^o's formula along flow of probability measures that we briefly recall (see \cite{cardelbook} for more details). Firstly, we fix some notations. Given $\mu\in\mathscr P_{\text{\tiny$2$}}(\R^n)$, for any $r\in[1,\infty)$ and $m\in\N\backslash\{0\}$, we use the shorthand notation $L_\mu^r(\R^m)$ to denote the space $L^r(\R^n,\Bc(\R^n),\mu;\R^m)$ of $\R^m$-valued $r$-integrable functions with respect to $\mu$. Similarly, $L_\mu^\infty(\R^m)$ denotes the space $L^\infty(\R^n,\Bc(\R^n),\mu;\R^m)$ of $\R^m$-valued $\mu$-essentially bounded functions.

Let $\vartheta$ be a real-valued function defined on  $\mathscr P_{\text{\tiny$2$}}(\R^n)$. Denote by $\upsilon$ the lifted version of $\vartheta$, that is the function defined on  $L^2(\Omega,\Gc,\P;\R^n)$ by $\upsilon(\xi)$ $=$ 
$\vartheta(\P_{\text{\tiny$\xi$}})$. We say that $\vartheta$ is differentiable (resp. $\Cc^1$) on
$\mathscr P_{\text{\tiny$2$}}(\R^n)$ if the lift $\upsilon$ is Fr\'echet differentiable (resp. continuously Fr\'echet differentiable) on  $L^2(\Omega,\Gc,\P;\R^n)$. 
In this case, the Fr\'echet derivative $[D \upsilon](\xi)$ of $\upsilon$ at $\xi$ $\in$ $L^2(\Omega,\Gc,\P;\R^n)$, viewed as an element $D\upsilon(\xi)$ of $L^2(\Omega,\Gc,\P;\R^n)$  by the Riesz representation theorem:  
$[D \upsilon](\xi)(Y)$ $=$ $\E[D\upsilon(\xi).Y]$ (we denote by $.$ the scalar product on $\R^n$),  can be represented as 
\beq \label{Uu1}
D\upsilon(\xi) &=& \partial_\mu \vartheta(\P_{_\xi})(\xi),
\enq
for some function  $\partial_\mu \vartheta(\mu)$ $:$ $\R^n$ $\rightarrow$ $\R^n$, with $\partial_\mu \vartheta(\mu)$ $\in$ $L^2_\mu(\R^n)$, depending only on the law  $\mu$ $=$ $\P_{_\xi}$ of $\xi$,   and  
called derivative of $\vartheta$ at $\mu$. We say that $\vartheta$ is partially  $\Cc^2$ if it is $\Cc^1$,  and one can find, for any $\mu$ $\in$ $\mathscr P_{\text{\tiny$2$}}(\R^n)$, 
a continuous version of the mapping $x\in\R^n$ $\mapsto$ $\partial_\mu \vartheta(\mu)(x)$, such that the mapping  
$(\mu,x)$ $\in$ $\mathscr P_{\text{\tiny$2$}}(\R^n)\times\R^n$ $\mapsto$ $\partial_\mu \vartheta(\mu)(x)$  is continuous at any point $(\mu,x)$ such that $x$ $\in$ Supp$(\mu)$, and if 
for each fixed $\mu$ $\in$ $\mathscr P_{\text{\tiny$2$}}(\R^n)$, the mapping $x$ $\in$ $\R^n$ $\mapsto$  $\partial_\mu \vartheta(\mu)(x)$ 
is differentiable in the standard sense, with a gradient denoted by   $\partial_x  \partial_\mu \vartheta(\mu)(x)$  $\in$ $\R^{n\times n}$, this derivative being jointly continuous at 
any $(\mu,x)$ $\in$ $\mathscr P_{\text{\tiny$2$}}(\R^n)\times\R^n$ such that $x$ $\in$ Supp$(\mu)$. 
We say that $\vartheta$ $\in$ $\Cc^2_b(\mathscr P_{\text{\tiny$2$}}(\R^n))$ if it is partially $\Cc^2$, $\partial_x  \partial_\mu \vartheta(\mu)$ $\in$ $L_\mu^\infty(\R^{n\times n})$, and for any compact set $\Kc$ of $\mathscr P_{\text{\tiny$2$}}(\R^n)$ we have
\beqs \label{rel1}
 \sup_{ \mu \in \Kc } \bigg[ \int_{\R^n} \big| \partial_\mu \vartheta(\mu)(x)\big|^2\mu(dx)  + 
\big\| \partial_x \partial_\mu \vartheta(\mu)\big\|_{_\infty}  
\bigg]  & < & \infty.
\enqs
Moreover, when the lifted function $\upsilon$ is twice continuously Fr\'echet differentiable, its second Fr\'echet derivative $D^2\upsilon(\xi)$, identified by the Riesz representation theorem as a self-adjoint (hence bounded) operator on $L^2(\Omega,\Gc,\P;\R^n)$, that is
$D^2\upsilon(\xi)$ $\in$ $S(L^2(\Omega,\Gc,\P;\R^n))$, is given by
\beq \label{rel2}
\E\big[ D^2\upsilon(\xi)(ZN).ZN \big] &=& \E \big[ {\rm tr}\big(\partial_x\partial_\mu \vartheta(\P_{\text{\tiny$\xi$}})(\xi)ZZ\trans \big)\big]  
\enq
for every $Z$ $\in$ $L^2(\Omega,\Gc,\P;\R^{n\times d})$ and any random vector $N$ $\in$ $L^2(\Omega,\Gc,\P;\R^d)$, with zero mean and unit variance, independent of $(\xi,Z)$.

Finally, we say that a function $\varphi$ $\in$ $C_b^{1,2}([0,T]\times \mathscr P_{\text{\tiny$2$}}(\R^n))$ if $\varphi$ is continuous on $[0,T]\times \mathscr P_{\text{\tiny$2$}}(\R^n)$, for every $t\in[0,T]$ the map $\varphi(t,\cdot)$ belongs to $\Cc^2_b(\mathscr P_{\text{\tiny$2$}}(\R^n))$, and for every $\mu\in\mathscr P_{\text{\tiny$2$}}(\R^n)$ the map $\varphi(\cdot,\mu)$ is continuously differentiable on $[0,T]$. Moreover, we say that a function $\phi$ $\in$ $C^{1,2}([0,T]\times L^2(\Omega,\Gc,\P;\R^n))$ if $\phi$ is continuous on $[0,T]\times L^2(\Omega,\Gc,\P;\R^n)$, for every $\xi\in L^2(\Omega,\Gc,\P;\R^n)$ the map $\phi(\cdot,\xi)$ is continuously Fr\'echet differentiable, for every $t\in[0,T]$ the map $\phi(t,\cdot)$ is twice continuously Fr\'echet differentiable.
 
\begin{Remark}
{\rm  The above definition of differentiability in the Wasserstein space, which is extrinsic via  the lifted identification with the Hilbert space of square-integrable random variables, is 
due to P.L. Lions \cite{lio12} (see also \cite{car13}),  
and turns out to be equivalent to a more intrinsic notion (of derivative in the Wasserstein space) used by various authors in connection with optimal transport  
and gradient flows (see e.g.  \cite{ambetal05}, \cite{ganswi15}), as recently shown in \cite{gantud17}.   
}
\epR
\end{Remark}

We also recall It\^o's formula along flow of probability measures for an It\^o process
\beqs
dX_t &=& b_t\,dt + \sigma_t\,dW_t, 
\enqs
where $b$ and $\sigma$ are $(\Fc_t\vee\Gc)_t$-progressively measurable processes. Then, for $\vartheta$ $\in$  $\Cc^2_b(\mathscr P_{\text{\tiny$2$}}(\R^n))$, we have
\beqs
\frac{d}{dt} \vartheta (\P_{_{X_t}}) &=& \E \bigg[ b_t.\partial_\mu \vartheta(\P_{_{X_t}})(X_t) + \frac{1}{2}{\rm tr}\big(\sigma_t\sigma_t\trans \partial_x\partial_\mu \vartheta(\P_{_{X_t}})(X_t) \big) \bigg]. 
\enqs

Let us now introduce the function defined on $\R^n\times\mathscr P_{\text{\tiny$2$}}(\R^n)\times A\times B\times\mathscr P(A\times B)\times\R^n\times\R^{n\times n}$  by 
\beqs
H(x,\mu,a,b,\nu,p,M) &=& \gamma(x,\mu,a,b,\nu).p + \frac{1}{2}{\rm tr}\big(\sigma\sigma\trans(x,\mu,a,b,\nu)M\big) + f(x,\mu,a,b,\nu).
\enqs

Finally, we denote by $L^0(\Omega,\Gc,\P;A)$ (resp. $L^0(\Omega,\Gc,\P;B)$) the set of $\Gc$-measurable random variables taking values in $A$ (resp. $B$).

\begin{Lemma}
Suppose that Assumption {\bf (A1)} holds and fix $\mu\in\mathscr P_{\text{\tiny$2$}}(\R^n)$. Consider two maps $\p\colon\R^n\rightarrow\R^n$ and $\M\colon\R^n\rightarrow\R^{n\times n}$ such that $\p\in L_\mu^2(\R^n)$ and $\M\in L_\mu^\infty(\R^{n\times n})$, then
\begin{align}\label{H=tildeH1}
\Sup_{\a \in L^0(\Omega,\Gc,\P;A)} \Inf_{\b \in L^0(\Omega,\Gc,\P;B)} 
\E \big[ H\big(\xi,\P_{_{\xi}},\a,\b,\P_{_{(\a,\b)}},\p(\xi),\M(\xi)\big) \big]& \notag \\
= \ \Sup_{\a \in L^0(\Omega,\Gc,\P;A)} \Inf_{\b \in L^0(\Omega,\Gc,\P;B)} 
\E \big[ H\big(\tilde\xi,\P_{_{\tilde\xi}},\a,\b,\P_{_{(\a,\b)}},\p(\tilde\xi),\M(\tilde\xi)\big) \big]&
\end{align}
and
\begin{align}\label{H=tildeH2}
\Inf_{\b \in L^0(\Omega,\Gc,\P;B)}   \Sup_{\a \in L^0(\Omega,\Gc,\P;A)} 
\E \big[ H\big(\xi,\P_{_{\xi}},\a,\b,\P_{_{(\a,\b)}},\p(\xi),\M(\xi)\big) \big]& \notag \\
= \ \Inf_{\b \in L^0(\Omega,\Gc,\P;B)}  \Sup_{\a \in L^0(\Omega,\Gc,\P;A)} 
\E \big[ H\big(\tilde\xi,\P_{_{\tilde\xi}},\a,\b,\P_{_{(\a,\b)}},\p(\tilde\xi),\M(\tilde\xi)\big) \big]&, 
\end{align}
for any $\xi,\tilde\xi\in L^2(\Omega,\Gc,\P;\R^n)$, with $\P_{\text{\tiny$\xi$}}=\P_{\text{\tiny$\tilde\xi$}}$. 
\end{Lemma}
{\bf Proof.}
We only prove \eqref{H=tildeH1}, the proof of \eqref{H=tildeH2} being analogous. In order to show \eqref{H=tildeH1}, we can proceed along the same lines as in the proof of Proposition \ref{P:Lift}, even though in this case the proof turns out to be much simpler (here $\a$ and $\b$ are random variables rather than controls or strategies as in Proposition \ref{P:Lift}). For this reason, we  skip the details and 
only report the main steps of the proof.

Similarly to the proof of Proposition \ref{P:Lift}, our aim is to prove the following:
\begin{align}\label{H_Htilde}
&\text{\emph{Given any $\a \in L^0(\Omega,\Gc,\P;A)$ and $\tilde\b \in L^0(\Omega,\Gc,\P;B)$, there exist $\tilde\a \in L^0(\Omega,\Gc,\P;A)$}} \notag \\
&\;\;\text{\emph{and $\b \in L^0(\Omega,\Gc,\P;B)$, with $\tilde\a$ (resp. $\b$) independent of $\tilde\b$ (resp. $\a$), such that:}} \notag \\
&\qquad\qquad\qquad\text{\emph{$(\xi,\a,\b)$ has the same law as }}(\tilde\xi,\tilde\a,\tilde\b)\text{\emph{, so that }} \\
&\qquad \E\big[ H\big(\xi,\P_{_{\xi}},\a,\b,\P_{_{(\a,\b)}},\p(\xi),\M(\xi)\big) \big] \ = \ \E\big[ H\big(\tilde\xi,\P_{_{\tilde\xi}},\tilde\a,\tilde\b,\P_{_{(\tilde\a,\tilde\b)}},\p(\tilde\xi),\M(\tilde\xi)\big) \big]. \notag
\end{align}
In order to prove \eqref{H_Htilde}, we recall that, by assumption, there exists a random variable $\Gamma^\Gc\colon(\Omega,\Gc)\rightarrow(G,\mathscr G)$, taking values in some Polish space $G$ with Borel $\sigma$-algebra $\mathscr G$, such that $\Gamma^\Gc$ has an atomless distribution and $\Gc=\sigma(\Gamma^\Gc)$. Since all atomless Polish probability spaces are isomorphic, we can suppose that $(G,\mathscr G)$ is $([0,1],\Bc([0,1]))$ and that $\Gamma^\Gc$ has uniform distribution.

Recalling that $\xi$ and $\a$ are $\Gc$-measurable, we have that
\[
\xi \ = \ \chi(\Gamma^\Gc), \qquad\qquad \a \ = \ \mathfrak a(\Gamma^\Gc), \qquad\qquad \P\text{-a.s.}
\]
for some measurable maps $\chi\colon[0,1]\rightarrow\R^n$ and $\mathfrak a\colon[0,1]\rightarrow A$. We can suppose, without loss of generality, that $\chi$ is surjective. Then, we notice that \eqref{H_Htilde} follows if we prove the following:
\begin{align}\label{tildeGamma^G}
\text{\emph{$\exists$ a $\Gc$-measurable random variable $\tilde\Gamma^\Gc$ with values in $[0,1]$ such that:}}& \\
\text{\emph{$\tilde\Gamma^\Gc$ has the same distribution of $\Gamma^\Gc$, $\sigma(\tilde\Gamma^\Gc)\vee\Nc=\Gc\vee\Nc$, and }} \tilde\xi \ = \ \chi(\tilde\Gamma^\Gc),\,\P\text{\emph{-a.s.}}& \notag
\end{align}
As a matter of fact, suppose that \eqref{tildeGamma^G} holds. Since $\sigma(\tilde\Gamma^\Gc)\vee\Nc=\Gc\vee\Nc$ and $\tilde\b$ is $\Gc$-measurable, we have
\[
\tilde\b \ = \ \tilde{\mathfrak b}(\tilde\Gamma^\Gc), \qquad \P\text{-a.s.}
\]
for some measurable function $\tilde{\mathfrak b}\colon[0,1]\rightarrow B$. Now, define
\[
\tilde\a \ := \ \mathfrak a(\tilde\Gamma^\Gc), \qquad\qquad \b \ := \ \tilde{\mathfrak b}(\Gamma^\Gc).
\]
Observe that $\tilde\a \in L^0(\Omega,\Gc,\P;A)$ and $\b \in L^0(\Omega,\Gc,\P;B)$. We also notice that $(\xi,\a,\b)$ has the same law as $(\tilde\xi,\tilde\a,\tilde\b)$, so that \eqref{H_Htilde} holds.

It remains to prove \eqref{tildeGamma^G}. Proceeding as in \textbf{Step IV} of the proof of Proposition \ref{P:Lift}, by the Jankov-von Neumann measurable selection theorem we deduce the existence of the analytically measurable right-inverse $\zeta$ of $\chi$. Then, we define
\[
\tilde\Gamma^\Gc \ := \ \zeta(\tilde\xi).
\]
Proceeding as in \textbf{Step IV} of the proof of Proposition \ref{P:Lift} we see that $\tilde\Gamma^\Gc$ satisfies \eqref{tildeGamma^G}, from which the claim follows.
\ep

\vspace{5mm}

In view of the above Lemma,  we can  define the lower and upper Hamiltonian functions $\mathscr H_-$, $\mathscr H_+$ $:$ $\mathscr P_{\text{\tiny$2$}}(\R^n)\times L_\mu^2(\R^n)\times L_\mu^\infty(\R^{n\times n})$ 
$\rightarrow$ $\R$ by 
\beqs
\mathscr H_-(\mu,\p,\M) &=&  \Sup_{\a \in L^0(\Omega,\Gc,\P;A)} \Inf_{\b \in L^0(\Omega,\Gc,\P;B)} 
\E \big[ H\big(\xi,\mu,\a,\b,\P_{_{(\a,\b)}},\p(\xi),\M(\xi)\big) \big]  \\
\mathscr H_+(\mu,\p,\M) &=&  \Inf_{\b \in L^0(\Omega,\Gc,\P;B)} \Sup_{\a \in L^0(\Omega,\Gc,\P;A)}  
\E \big[ H\big(\xi,\mu,\a,\b,\P_{_{(\a,\b)}},\p(\xi),\M(\xi)\big) \big]    
\enqs
for every $(\mu,\p,\M)$ $\in$ $\mathscr P_{\text{\tiny$2$}}(\R^n)\times L_\mu^2(\R^n)\times L_\mu^\infty(\R^{n\times n})$, with  $\xi$ $\in$ 
$L^2(\Omega,\Gc,\P;\R^n)$ such that $\P_{_{\xi}}$ $=$ $\mu$. Then, consider the lower Bellman-Isaacs equation on $[0,T]\times\mathscr P_{\text{\tiny$2$}}(\R^n)$:
\begin{equation}\label{lowerBI}
\begin{cases}
\vspace{1mm}\displaystyle- \Dt{\vartheta}(t,\mu) -  \mathscr H_-(\mu,\partial_\mu\vartheta(t,\mu),\partial_x\partial_\mu\vartheta(t,\mu)) &= \ 0,  \hspace{2.8cm}  (t,\mu) \in [0,T) \times \mathscr P_{\text{\tiny$2$}}(\R^n), \\
\hspace{5cm}\vartheta(T,\mu) &= \ \displaystyle \int_{\R^n} g(x,\mu) \mu(dx), \;\;\;\;\;\!\! \mu \in  \mathscr P_{\text{\tiny$2$}}(\R^n)
\end{cases}
\end{equation}
and the upper Bellman-Isaacs equation on $[0,T]\times\mathscr P_{\text{\tiny$2$}}(\R^n)$:
\begin{equation}\label{upperBI}
\begin{cases}
\vspace{2mm}\displaystyle- \Dt{\vartheta}(t,\mu) -  \mathscr H_+(\mu,\partial_\mu\vartheta(t,\mu),\partial_x\partial_\mu\vartheta(t,\mu)) &= \ 0,  \hspace{2.8cm}  (t,\mu) \in [0,T) \times \mathscr P_{\text{\tiny$2$}}(\R^n), \\
\hspace{5cm}\vartheta(T,\mu) &= \ \displaystyle \int_{\R^n} g(x,\mu) \mu(dx), \;\;\;\;\;\!\! \mu \in  \mathscr P_{\text{\tiny$2$}}(\R^n).
\end{cases}
\end{equation}

The corresponding lifted equation on $[0,T]\times L^2(\Omega,\Gc,\P;\R^n)$ is formulated as follows in view of the relations \reff{Uu1}-\reff{rel2} 
between derivatives in the Wasserstein space and Fr\'echet derivatives. We define the functions $\Hc_-$, $\Hc_+$ $:$ $L^2(\Omega,\Gc,\P;\R^n)\times L^2(\Omega,\Gc,\P;\R^n)\times  S(L^2(\Omega,\Gc,\P;\R^n))$  $\rightarrow$ $\R$ by
\beqs
\Hc_-(\xi,P,Q) &=& \Sup_{\a \in L^0(\Omega,\Gc,\P;A)} \Inf_{\b \in L^0(\Omega,\Gc,\P;B)} 
\E\Big[  \gamma(\xi,\P_{\text{\tiny$\xi$}},\a,\b,\P_{_{(\a,\b)}}).P  \\
& & \hspace{1mm} + \; \frac{1}{2} Q(\sigma(\xi,\P_{\text{\tiny$\xi$}},\a,\b,\P_{_{(\a,\b)}}) N). \sigma(\xi,\P_{\text{\tiny$\xi$}},\a,\b,\P_{_{(\a,\b)}})N  
\; + \; f(\xi,\P_{\text{\tiny$\xi$}},\a,\b,\P_{_{(\a,\b)}}) \Big], \\
\Hc_+(\xi,P,Q) &=& \Inf_{\b \in L^0(\Omega,\Gc,\P;B)}  \Sup_{\a \in L^0(\Omega,\Gc,\P;A)}  
\E\Big[  \gamma(\xi,\P_{\text{\tiny$\xi$}},\a,\b,\P_{_{(\a,\b)}}).P  \\
& & \hspace{1mm} + \; \frac{1}{2} Q(\sigma(\xi,\P_{\text{\tiny$\xi$}},\a,\b,\P_{_{(\a,\b)}}) N). \sigma(\xi,\P_{\text{\tiny$\xi$}},\a,\b,\P_{_{(\a,\b)}})N  
\; + \; f(\xi,\P_{\text{\tiny$\xi$}},\a,\b,\P_{_{(\a,\b)}}) \Big] 
\enqs
where $N$ $\in$ $L^2(\Omega,\Gc,\P;\R^d)$, with zero mean and unit variance, is independent of $\xi$. Then, consider the lower Bellman-Isaacs equation on  $[0,T]\times L^2(\Omega,\Gc,\P;\R^n)$:
\begin{equation}\label{lowerBIlift}
\begin{cases}
\vspace{1mm}\displaystyle - \Dt{\upsilon}(t,\xi) -  \Hc_-(\xi,D\upsilon(t,\xi),D^2 \upsilon(t,\xi))\!\!\!\! &= \ 0,  \hspace{2cm}  (t,\xi) \in [0,T) \times L^2(\Omega,\Gc,\P;\R^n), \\
\hspace{4.4cm}\upsilon(T,\xi)\!\!\!\! &= \ \displaystyle \E\big[ g(\xi,\P_{\text{\tiny$\xi$}})\big], \;\;\;\;\; \xi \in  L^2(\Omega,\Gc,\P;\R^n)
\end{cases}
\end{equation}
and the upper Bellman-Isaacs equation on $[0,T]\times L^2(\Omega,\Gc,\P;\R^n)$:
\begin{equation}\label{upperBIlift}
\begin{cases}
\vspace{1mm}\displaystyle - \Dt{\upsilon}(t,\xi) -  \Hc_+(\xi,D\upsilon(t,\xi),D^2 \upsilon(t,\xi))\!\!\!\! &= \ 0,  \hspace{2cm}  (t,\xi) \in [0,T) \times L^2(\Omega,\Gc,\P;\R^n), \\
\hspace{4.4cm}\upsilon(T,\xi)\!\!\!\! &= \ \displaystyle \E\big[ g(\xi,\P_{\text{\tiny$\xi$}})\big], \;\;\;\;\; \xi \in  L^2(\Omega,\Gc,\P;\R^n).
\end{cases}
\end{equation}

\begin{Remark}
{\rm In the case where the coefficients do not depend on the law of the controls, i.e. $\gamma$ $=$ $\gamma(x,\mu,a,b)$, $\sigma$ $=$ $\sigma(x,\mu,a,b)$, $f$ $=$ 
$f(x,\mu,a,b)$, so that $H$ $=$ $H(x,\mu,a,b,p,M)$, the optimization over $\a$ $\in$ $L^0(\Omega,\Gc,\P;A)$ and $\b$ $\in$ 
$L^0(\Omega,\Gc,\P;B)$ in the Hamiltonian functions $\mathscr H_-$ and 
$\mathscr H_+$ reduces to a pointwise optimization over $A$ and $B$ inside the expectation operator, namely
\beq
\mathscr H_-(\mu,\p,\M) &=& \E \Big[  \Sup_{a\in A} \Inf_{b \in B} H(\xi,\mu,a,b,\p(\xi),\M(\xi)) \Big]  \label{HnoMF} \\
\mathscr H_+(\mu,\p,\M) &=& \E \Big[  \Inf_{b \in B}  \Sup_{a\in A}  H(x,\mu,a,b,\p(\xi),\M(\xi)) \Big].  \nonumber
\enq
Indeed, it is clear that for any $\a \in L^0(\Omega,\Gc,\P;A)$, 
\beqs
\Inf_{\b \in L^0(\Omega,\Gc,\P;B)} \E \big[ H\big(\xi,\mu,\a,\b,\p(\xi),\M(\xi)\big) \big]  & \geq &  
\E \big[ \inf_{b \in B} H\big(\xi,\mu,\a,b,\p(\xi),\M(\xi)\big) \big] 
\enqs
Conversely,  by the Jankov-von Neumann measurable selection theorem, for any $\eps$ $>$ $0$, there exists a measurable function $x$ $\in$ $\R^n$ $\mapsto$ 
$b^\eps(x)$ $\in$ $B$ (depending on ${\bf a}$) such that $\Inf_{b \in B} H(t,x,\mu,\a,b,\p(x),\M(x))$ $\geq$ $H(t,x,\mu,\a,b^\eps(x),\p(x),\M(x))$ $-$ $\eps$, $\mu(dx)$-a.e.. 
By considering ${\bf b}^\eps$ $=$ $b^\eps(\xi)$ $\in$ $L^0(\Omega,\Gc,\P;B)$, we then have 
\beqs
\E \big[ \inf_{b \in B} H\big(\xi,\mu,\a,b,\p(\xi),\M(\xi)\big) \big] & \geq & \E \big[ H\big(\xi,\mu,\a,\b^\eps,\p(\xi),\M(\xi)\big) \big]  - \eps \\
& \geq & \Inf_{\b \in L^0(\Omega,\Gc,\P;B)} \E \big[ H\big(\xi,\mu,\a,\b,\p(\xi),\M(\xi)\big) \big]  - \eps, 
\enqs
and thus, by sending $\eps$ to zero, the equality: 
\beqs
\Inf_{\b \in L^0(\Omega,\Gc,\P;B)} \E \big[ H\big(\xi,\mu,\a,\b,\p(\xi),\M(\xi)\big) \big]  & = &  
\E \big[ \inf_{b \in B} H\big(\xi,\mu,\a,b,\p(\xi),\M(\xi)\big) \big].
\enqs 
Next, following the same argument as above, we show that
\beqs
 \Sup_{\a \in L^0(\Omega,\Gc,\P;A)} \Inf_{\b \in L^0(\Omega,\Gc,\P;B)} \E \big[ H\big(\xi,\mu,\a,\b,\p(\xi),\M(\xi)\big) \big]  & = &  
\E \big[ \sup_{a\in A} \inf_{b \in B} H\big(\xi,\mu,a,b,\p(\xi),\M(\xi)\big) \big] 
\enqs
Similarly, the optimization over $\a$ $\in$ $L^0(\Omega,\Gc,\P;A)$ and $\b$ $\in$ 
$L^0(\Omega,\Gc,\P;B)$ in the Hamiltonian functions $\Hc_-$ and 
$\Hc_+$ reduces to a pointwise optimization over $A$ and $B$ inside the expectation operator.  
}
\epR
\end{Remark}

\vspace{3mm}

\begin{Remark} \label{remHJBI} 
{\rm (Case without mean-field interaction) 

\noindent Consider the standard zero-sum stochastic differential game as in \cite{FlemingSouganidis} where all the coefficients depend only on state and controls. 
Given a function $v^{FS}$ $\in$ $C^{1,2}([0,T]\times\R^n)$ with quadratic growth condition on its derivatives,  let us  define the function 
$\vartheta$ on $[0,T]\times \mathscr P_{\text{\tiny$2$}}(\R^n)$ by 
\beqs
\vartheta(t,\mu) &=& \E[ v^{FS}(t,\xi) ], 
\;\;\;\quad \mu \in  \mathscr P_{\text{\tiny$2$}}(\R^n), \; \xi \in L^2(\Omega,\Gc,\P;\R^n), \; \P_{_{\xi}} = \mu. 
\enqs
Then $\vartheta(t,\cdot)$ is partially $\Cc^2$, while $\vartheta(\cdot,\mu)$ is $C^1$, with
\[
\partial_\mu \vartheta(t,\mu) \ = \ D_x v^{FS}(t,\cdot), \qquad \partial_x \partial_\mu \vartheta(t,\mu) \ = \ D_x^2 v^{FS}(t,\cdot), \qquad \Dt{\vartheta}(t,\mu) \ = \ \E\bigg[\Dt{v^{FS}}(t,\xi)\bigg].
\]
Moreover, we have (recalling also \reff{HnoMF})
\begin{align} 
&\Dt{\vartheta}(t,\mu) +  \mathscr H_-(\mu,\partial_\mu\vartheta(t,\mu),\partial_x\partial_\mu\vartheta(t,\mu))  \nonumber \\
&= \ \E \bigg[    \Dt{v^{FS}}(t,\xi)  +  H_-^{FS}\big(\xi,D_x v^{FS}(t,\xi),D_x^2 v^{FS}(t,\xi) \big) \bigg], \label{connecHJBI}
\end{align}
where $H_-^{FS}(x,p,M)$ $=$ $\sup_{a \in A} \inf_{b \in B}\big[ \gamma(x,a,b).p + \frac{1}{2}{\rm tr}(\sigma\sigma\trans(x,a)M \big) + f(x,a,b) \big]$ is the lower Bellman-Isaacs Hamiltonian associated to 
the zero-sum stochastic differential game as in \cite{FlemingSouganidis}.  The connection \reff{connecHJBI}  
between the lower Bellman-Isaacs equation on $[0,T]\times\R^n$ and the lower Bellman-Isaacs equation on  $[0,T]\times  \mathscr P_{\text{\tiny$2$}}(\R^n)$ shows  that 
$v^{FS}$ is a  (smooth) solution to 
\[
\begin{cases}
\vspace{2mm}
\displaystyle - \Dt{v^{FS}}(t,x)  -  H_-^{FS}\big(x,D_x v^{FS}(t,x),D_x^2 v^{FS}(t,x) \big) = 0, &\qquad (t,x) \in [0,T)\times\R^n, \\
v^{FS}(T,x) = g(x), &\qquad x \in \R^n, 
\end{cases}
\]
if and only if $\vartheta$  is a (smooth) solution to \reff{lowerBI}. A similar  connection holds for the upper Bellman-Isaacs equation. 
}
\epR
\end{Remark}

\vspace{5mm}

We now consider two definitions of viscosity solution for the Bellman-Isaacs equations, on one hand on the Wasserstein space  $\mathscr P_{\text{\tiny$2$}}(\R^n)$ and on the other hand on the lifted Hilbert space $L^2(\Omega,\Gc,\P;\R^n)$. In the sequel, the Hamiltonian function $\mathscr H$ denotes either $\mathscr H_-$ or $\mathscr H_+$, and $\Hc$ stands for $\Hc_-$ or $\Hc_+$.

\begin{Definition}[Viscosity solution in $\mathscr P_{\text{\tiny$2$}}$]

\noindent  A continuous function $\vartheta$ on $[0,T]\times\mathscr P_{\text{\tiny$2$}}(\R^n)$ is a viscosity solution to \reff{lowerBI} (or \reff{upperBI}) if:
\begin{itemize}
\item[\textup{(i)}] \textup{(viscosity supersolution property):} $\vartheta(T,\mu)$ $\geq$ $\int_{\R^n} g(x,\mu) \mu(dx)$,  $\mu$ $\in$ $\mathscr P_{\text{\tiny$2$}}(\R^n)$, and for any test function 
$\varphi$ $\in$ $C_b^{1,2}([0,T]\times \mathscr P_{\text{\tiny$2$}}(\R^n))$ such that $\vartheta - \varphi$ has a minimum at $(t_{_0},\mu_{_0})$ $\in$ $[0,T)\times\mathscr P_{\text{\tiny$2$}}(\R^n)$, 
we have 
\beqs
- \Dt{\varphi}(t_{_0},\mu_{_0}) -  \mathscr H(\mu_{_0},\partial_\mu\varphi(t_{_0},\mu_{_0}),\partial_x\partial_\mu\varphi(t_{_0},\mu_{_0})) & \geq & 0. 
\enqs 
\item[\textup{(ii)}] \textup{(viscosity subsolution property):} $\vartheta(T,\mu)$ $\leq$ $\int_{\R^n} g(x,\mu) \mu(dx)$,  $\mu$ $\in$ $\mathscr P_{\text{\tiny$2$}}(\R^n)$, and for any test function 
$\varphi$ $\in$ $C_b^{1,2}([0,T]\times \mathscr P_{\text{\tiny$2$}}(\R^n))$ such that $\vartheta - \varphi$ has a maximum at $(t_{_0},\mu_{_0})$ $\in$ $[0,T)\times\mathscr P_{\text{\tiny$2$}}(\R^n)$, we have 
\beqs
- \Dt{\varphi}(t_{_0},\mu_{_0}) -  \mathscr H(\mu_{_0},\partial_\mu\varphi(t_{_0},\mu_{_0}),\partial_x\partial_\mu\varphi(t_{_0},\mu_{_0})) & \leq & 0.
\enqs
\end{itemize}
\end{Definition}

\vspace{5mm}

\begin{Definition}[Viscosity solution in $L^2$]\label{defviscoL2}

\noindent  A continuous function $\upsilon$ on $[0,T]\times L^2(\Omega,\Gc,\P;\R^n)$ is a viscosity solution to \reff{lowerBIlift} (or \reff{upperBIlift}) if:
\begin{itemize}
\item[\textup{(i)}] \textup{(viscosity supersolution property):} $\upsilon(T,\xi)$ $\geq$ $\E[ g(\xi,\P_{\text{\tiny$\xi$}})]$,  $\xi$ $\in$ $L^2(\Omega,\Gc,\P;\R^n)$, and for any test function 
$\phi$ $\in$ $C^{1,2}([0,T]\times L^2(\Omega,\Gc,\P;\R^n))$ such that $\upsilon - \phi$ has a minimum at $(t_{_0},\xi_{_0})$ $\in$ $[0,T)\times L^2(\Omega,\Gc,\P;\R^n)$, we have 
\beqs
- \Dt{\phi}(t_{_0},\xi_{_0}) -  \Hc(\xi_{_0},D\phi(t_{_0},\xi_{_0}),D^2 \phi(t_{_0},\xi_{_0})) & \geq & 0.
\enqs
\item[\textup{(ii)}] \textup{(viscosity subsolution property):} $\upsilon(T,\xi)$ $\leq$ $\E[ g(\xi,\P_{\text{\tiny$\xi$}})]$,  $\xi$ $\in$ $L^2(\Omega,\Gc,\P;\R^n)$, and for any test function 
$\phi$ $\in$ $C^{1,2}([0,T]\times L^2(\Omega,\Gc,\P;\R^n))$ such that $\upsilon - \phi$ has a maximum at $(t_{_0},\xi_{_0})$ $\in$ $[0,T)\times L^2(\Omega,\Gc,\P;\R^n)$, we have 
\beqs
- \Dt{\phi}(t_{_0},\xi_{_0}) -  \Hc(\xi_{_0},D\phi(t_{_0},\xi_{_0}),D^2 \phi(t_{_0},\xi_{_0})) & \leq & 0.
\enqs
\end{itemize}
\end{Definition}

\begin{Remark} \label{remvisco}
{\rm Given these two definitions of 
viscosity solutions in $\mathscr P_{\text{\tiny$2$}}(\R^n)$ and in $L^2(\Omega,\Gc,\P;\R^n)$, a natural question is the connection between the viscosity property of 
$\vartheta$ to the Bellman-Isaacs equation  \reff{lowerBI} (or \reff{upperBI}) and the viscosity property of its lifted function $\upsilon$ to \reff{lowerBIlift} (or \reff{upperBIlift}). 
Actually, as pointed out in \cite{bucetal17} (see their Example 2.1), for a  function $\varphi$ $\in$ $C_b^{1,2}([0,T]\times \mathscr P_{\text{\tiny$2$}}(\R^n))$ its lifted function $\phi$ may not be in general in $C^{1,2}([0,T]\times L^2(\Omega,\Gc,\P;\R^n))$, which means that we cannot deduce from the viscosity property of $\upsilon$ to \reff{lowerBIlift} (or \reff{upperBIlift}) the viscosity 
property of $\vartheta$ to  \reff{lowerBI} (or \reff{upperBI}).  On the other hand, since a test function $\phi$ $\in$ $C^{1,2}\big([0,T]\times L^2(\Omega,\Gc,\P;\R^n)\big)$ is in general not necessarily the lifted function of a test function $\varphi$ in  $[0,T]\times\mathscr P_{\text{\tiny$2$}}(\R^n)$, we cannot claim that the viscosity property of $\vartheta$ to  \reff{lowerBI} (or \reff{upperBI}) implies the viscosity property of its lifted function $\upsilon$ to \reff{lowerBIlift} (or \reff{upperBIlift}).  This would hold true whenever we restrict in Definition \ref{defviscoL2} 
to test functions $\phi$ $\in$  $C^{1,2}([0,T]\times L^2(\Omega,\Gc,\P;\R^n))$ such that $\phi(t,\xi)$ depends on $\xi$ only via its law $\P_{_\xi}$, hence are lifted from functions in  
$C_b^{1,2}([0,T]\times \mathscr P_{\text{\tiny$2$}}(\R^n))$. However, in this case, we could not rely  on general comparison principles  for viscosity solutions in Hilbert spaces 
(see \cite{fabgozswi17}), which is in fact our main motivation for the introduction of Definition \ref{defviscoL2}, see Remark \ref{remunivisco}. 
}
\epR
\end{Remark}

From the DPP, we can now prove the viscosity solution property of the lower and upper value functions to the lower and upper Bellman-Isaacs equations. 

\begin{Theorem} \label{prop:visco}
Let Assumptions {\bf (A1)} and {\bf (A2)} hold. 

\noindent\textup{1)} The inverse-lifted  lower (resp. upper) value function $\mathscr V$ (resp. $\mathscr U$) is a viscosity solution to the lower (resp. upper) Bellman-Isaacs equation \reff{lowerBI} (resp. \reff{upperBI}) on  
$[0,T]\times\mathscr P_{\text{\tiny$2$}}(\R^n)$.

\noindent\textup{2)} The lower (resp. upper) value function $v$ (resp. $u$) is a viscosity solution to the lower (resp. upper) Bellman-Isaacs equation 
\reff{lowerBIlift} (resp. \reff{upperBIlift}) on  $[0,T]\times L^2(\Omega,\Gc,\P;\R^n)$.
\end{Theorem}
\noindent {\bf Proof.}
We only prove result 1), that is the viscosity property in $\mathscr P_{\text{\tiny$2$}}(\R^n)$ (as the viscosity property in $L^2(\Omega,\Gc,\P;\R^n)$ has a similar proof), and for the  
inverse-lifted lower value function $\mathscr V$, as the proof for the inverse-lifted upper value function is analogous. Obviously, $\mathscr V(T,\mu)$ $=$ $\E[g(\xi,\P_{_\xi})]$ $=$ $\int_{\R^n} g(x,\mu) \mu(dx)$, for $\mu$ $\in$  $\mathscr P_{\text{\tiny$2$}}(\R^n)$ and   
$\xi$ $\in$ $L^2(\Omega,\Gc,\P;\R^n)$ such that $\P_{_\xi}$ $=$ $\mu$. 

\vspace{2mm}

\noindent (i) {\it Viscosity supersolution property}.  
Fix $(t_{_0},\mu_{_0})$ $\in$ $[0,T)\times\mathscr P_{\text{\tiny$2$}}(\R^n)$, $\xi_{_0}$ $\in$ $L^2(\Omega,\Gc,\P;\R^n)$ with $\P_{_{\xi_{_0}}}$ $=$ $\mu_{_0}$, 
and consider any test function 
$\varphi$ $\in$ $C_b^{1,2}([0,T]\times \mathscr P_{\text{\tiny$2$}}(\R^n))$ such that $\mathscr V - \varphi$ has a minimum at $(t_{_0},\mu_{_0})$, and w.l.o.g. 
$\mathscr V(t_{_0},\mu_{_0})$ $=$   $\varphi(t_{_0},\mu_{_0})$.  We argue by contradiction, and assume on the contrary that 
\beqs
\Sup_{\a \in L^0(\Omega,\Gc,\P;A)} \Inf_{\b \in L^0(\Omega,\Gc,\P;B)}  \E \big[ F_\varphi(t_{_0},\xi_{_0},\mu_{_0},\a,\b,\P_{(\a,\b)}) \big] & \geq & 3\,\eps, 
\enqs
for some $\eps$ $>$ $0$, where we set 
\beq \label{defF}
F_\varphi(t,x,\mu,a,b,\nu) &=& \Dt{\varphi}(t,\mu) + H\big(x,\mu,a,b,\nu,\partial_\mu \varphi(t,\mu)(x),\partial_x\partial_\mu \varphi(t,\mu)(x) \big).   
\enq
This implies that there exists $\a^\eps$ $\in$ $L^0(\Omega,\Gc,\P;A)$ such that for all $\b \in L^0(\Omega,\Gc,\P;B)$
\beq \label{Fepssur}
\E \big[ F_\varphi(t_{_0},\xi_{_0},\mu_{_0},\a^\eps,\b,\P_{(\a^\eps,\b)}) \big] & \geq & 2\,\eps. 
\enq
By Remark \ref{R:ExH}, it follows that \eqref{Fepssur2} holds for any $\b \in L^0(\Omega,\Fc,\P;B)$. As a matter of fact, take $\b \in L^0(\Omega,\Fc,\P;B)$ and denote by $\pi$ the distribution on $(\R^n\times A\times B,\Bc(\R^n)\otimes\Bc(A)\otimes\Bc(B))$ of the random vector $(\xi_{_0},\a^\eps,\b)$. From Remark \ref{R:ExH} (with $E$, $H$, $\zeta$ corresponding respectively to $\R^n\times A$, $B$, $(\xi_{_0},\a^\eps)$) we deduce the existence of a measurable map $\bar\b\colon(\Omega,\Gc)\rightarrow(B,\Bc(B))$ (the map $\bar\b$ corresponds to the map $\eta$ in Remark \ref{R:ExH}) such that $\P_{\text{\tiny$(\xi_{_0},\a^\eps,\bar\b)$}}=\pi$. So, in particular,
\[
\E \big[ F_\varphi(t_{_0},\xi_{_0},\mu_{_0},\a^\eps,\b,\P_{(\a^\eps,\b)}) \big] \ = \ \E \big[ F_\varphi(t_{_0},\xi_{_0},\mu_{_0},\a^\eps,\bar\b,\P_{(\a^\eps,\bar\b)}) \big] \ \geq \ 2\,\eps.
\]
Hence, \eqref{Fepssur2} holds for any $\b \in L^0(\Omega,\Fc,\P;B)$. Now, under the continuity assumptions {\bf (A1)} and {\bf (A2)}, we easily see that $H(x,\mu,a,b,\nu,p,M)$ is continuous in $(x,\mu,p,M)$ uniformly with respect to $(a,b,\nu)$, from which we deduce,  recalling the continuity of $\partial_\mu\varphi(t,\mu)(x)$ and $\partial_x\partial_\mu\varphi(t,\mu)(x)$, and by Lebesgue's dominated convergence theorem, 
the continuity of $(t,\xi,\mu)$ $\in$ $[0,T]\times L^2(\Omega,\Fc_T\vee\Gc,\P;\R^n)\times\mathscr P_{\text{\tiny$2$}}(\R^n)$ $\mapsto$ 
$\E\big[ F_\varphi(t,\xi,\mu,\a,\b,\P_{(\a,\b)}) \big]$ uniformly with respect to $\a$ $\in$ $\in$ $L^0(\Omega,\Gc,\P;A)$ and $\b$ $\in$ $L^0(\Omega,\Fc,\P;B)$.  From \reff{Fepssur},  there is some $\delta$ $>$ $0$ such that for all $\b \in L^0(\Omega,\Fc,\P;B)$
\beq \label{Fepssur2}
\E \big[ F_\varphi(s,\xi,\mu,\a^\eps,\b,\P_{(\a^\eps,\b)}) \big] & \geq & \eps, \qquad \forall s \in [t_{_0},t_{_0}+\delta], \; (\xi,\mu) \in B_\delta(\xi_{_0},\mu_{_0}).
\enq
Here $B_\delta(\xi_{_0},\mu_{_0})$ is the ball of center $(\xi_{_0},\mu_{_0})$ and  radius $\delta$ in the metric space 
$L^2(\Omega,\Fc_T\vee\Gc,\P;\R^n)\times\mathscr P_{\text{\tiny$2$}}(\R^n)$.  
From \reff{estimWasser}-\reff{estimX1}, we have for all $s$ $\in$ $[t_{_0},T]$, $\alpha$ $\in$ $\Ac$, $\beta$ $\in$ $\Bc$, 
\beqs
\E\big| X_s^{t_{_0},\xi_{_0},\alpha,\beta} - \xi_{_0}\big|^2 + \Wc_{_2}^2\big(\P_{\text{\tiny$X_s^{t_{_0},\xi_{_0},\alpha,\beta}$}},\mu_{_0}\big) & \leq & C\,\big(1 + \E|\xi_{_0}|^2\big)\,|s-t_{_0}|.  
\enqs
We can then pick some $h$ $>$ $0$ small enough  so that for all $\alpha$ $\in$ $\Ac$, $\beta$ $\in$ $\Bc$
\beq \label{hpetit}
\big(s,X_s^{t_{_0},\xi_{_0},\alpha,\beta},\P_{\text{\tiny$X_s^{t_{_0},\xi_{_0},\alpha,\beta}$}} \big) \in [t_{_0},t_{_0}+\delta]\times B_\delta(\xi_{_0},\mu_{_0}), \;\;\; 
\mbox{ for all } s \in [t_{_0},t_{_0} + h]. 
\enq
Consider the constant control $\alpha^\eps$ in $\Ac$ equal to $\a^\eps$, and take an  arbitrary  $\beta[\cdot]\in\Bc_{\textup{\tiny str}}$. By applying It\^o's formula to 
$\varphi(r,\P_{\text{\tiny$X_{r}^{t_{_0},\xi_{_0},\alpha^\eps,\beta[\alpha^\eps]}$}})$ between  $t_{_0}$ and $t_{_0}+h$, we then get  by \reff{Fepssur2}-\reff{hpetit} 
\begin{align*}
&\varphi\big(t_{_0}+h,\P_{\text{\tiny$X_{t_{_0}+h}^{t_{_0},\xi_{_0},\alpha^\eps,\beta[\alpha^\eps]}$}}\big) -   \varphi(t_{_0},\mu_{_0}) \\
&= \ \E\bigg[\int_{t_{_0}}^{t_{_0}+h} 
\Big(F_\varphi\big(s,X_s^{t_{_0},\xi_{_0},\alpha^\eps,\beta[\alpha^\eps]},\P_{\text{\tiny$X_{s}^{t_{_0},\xi_{_0},\alpha^\eps,\beta[\alpha^\eps]}$}},
\a^\eps,\beta[\alpha^\eps]_s,\P_{\text{\tiny$(\a^\eps,\beta[\alpha^\eps]_s)$}} \big)  \\
&\quad \ - f\big(X_s^{t_{_0},\xi_{_0},\alpha^\eps,\beta[\alpha^\eps]},\P_{\text{\tiny$X_s^{t_{_0},\xi_{_0},\alpha^\eps,\beta[\alpha^\eps]}$}},\a^\eps,\beta[\alpha^\eps]_s,
\P_{\text{\tiny$(\a^\eps,\beta[\alpha^\eps]_s)$}}\big)\Big)\,ds \bigg]  \\
&\geq \ \eps h - \E \bigg[  \int_{t_{_0}}^{t_{_0}+h}  f\big(X_s^{t_{_0},\xi_{_0},\alpha^\eps,\beta[\alpha^\eps]},\P_{\text{\tiny$X_s^{t_{_0},\xi_{_0},\alpha^\eps,\beta[\alpha^\eps]}$}},\a^\eps,\beta[\alpha^\eps]_s,\P_{\text{\tiny$(\a^\eps,\beta[\alpha^\eps]_s)$}}\big)\,ds \bigg]. 
\end{align*}
Recalling that $\mathscr V(t_{_0},\mu_{_0})$ $=$   $\varphi(t_{_0},\mu_{_0})$, $\mathscr V$ $\geq$ $\varphi$, and that $\beta$ is arbitrary in $\Bc_{\textup{\tiny str}}$, we obtain
\begin{align*}
\mathscr V(t_{_0},\mu_{_0}) + \eps h \ \leq \  \inf_{\beta[\cdot]\in\Bc_{\textup{\tiny str}}}  \sup_{\alpha\in\Ac}  \E \bigg[  \int_{t_{_0}}^{t_{_0}+h}  f\big(X_s^{t_{_0},\xi_{_0},\alpha,\beta[\alpha]},\P_{\text{\tiny$X_s^{t_{_0},\xi_{_0},\alpha,\beta[\alpha]}$}},\alpha_s,\beta[\alpha]_s,\P_{\text{\tiny$(\alpha,\beta[\alpha]_s)$}}\big)\,ds \\
+ \mathscr V\big(t_{_0}+h,\P_{\text{\tiny$X_{t_{_0}+h}^{t_{_0},\xi_{_0},\alpha,\beta[\alpha]}$}}\big) \bigg], 
\end{align*}
which contradicts the DPP relation \reff{DPPvinv}. 
 
\vspace{1mm}

\noindent (ii) {\it Viscosity subsolution property}.  Fix $(t_{_0},\mu_{_0})$ $\in$ $[0,T)\times\mathscr P_{\text{\tiny$2$}}(\R^n)$, $\xi_{_0}$ $\in$ $L^2(\Omega,\Gc,\P;\R^n)$ with $\P_{_{\xi_{_0}}}$ $=$ $\mu_{_0}$, and consider any test function 
$\varphi$ $\in$ $C_b^{1,2}([0,T]\times \mathscr P_{\text{\tiny$2$}}(\R^n))$ such that $\mathscr V - \varphi$ has a maximum at $(t_{_0},\mu_{_0})$, and w.l.o.g. 
$\mathscr V(t_{_0},\mu_{_0})$ $=$   $\varphi(t_{_0},\mu_{_0})$. We still argue by contradiction, and assume on the contrary that 
\beq\label{psi1bis}
\Sup_{\a \in L^0(\Omega,\Gc,\P;A)} \Inf_{\b \in L^0(\Omega,\Gc,\P;B)}  \E \big[ F_\varphi(t_{_0},\xi_{_0},\mu_{_0},\a,\b,\P_{(\a,\b)}) \big] & \leq & - 3\,\eps, 
\enq
for some $\eps$ $>$ $0$, where $F_\varphi$ is defined as in \reff{defF}. As for \eqref{Fepssur}, by Remark \ref{R:ExH} it follows that \eqref{psi1bis} still holds if we take the supremum over all $\a \in L^0(\Omega,\Fc_T\vee\Gc,\P;A)$ (actually, over all $\a \in L^0(\Omega,\Fc,\P;A)$, but here we take $\Fc_T\vee\Gc$ since it is a countably generated $\sigma$-algebra). As a matter of fact, given $\a \in L^0(\Omega,\Fc_T\vee\Gc,\P;A)$ and $\b\in L^0(\Omega,\Gc,\P;B)$, denoting $\pi=\P_{\text{\tiny$(\xi_{_0},\a,\b)$}}$, by Remark \ref{R:ExH} there exists $\bar\a \in L^0(\Omega,\Gc,\P;A)$ such that $\P_{\text{\tiny$(\xi_{_0},\bar\a,\b)$}}=\pi$; so, in particular,
\[
\E \big[ F_\varphi(t_{_0},\xi_{_0},\mu_{_0},\a,\b,\P_{(\a,\b)}) \big] \ = \ \E \big[ F_\varphi(t_{_0},\xi_{_0},\mu_{_0},\bar\a,\b,\P_{(\bar\a,\b)}) \big].
\]
Taking the infimum over $\b\in L^0(\Omega,\Gc,\P;B)$, and then the supremum over $\a \in L^0(\Omega,\Fc_T\vee\Gc,\P;A)$, by \eqref{psi1bis} we end up with
\begin{equation}\label{psi1}
\Sup_{\a \in L^0(\Omega,\Fc_T\vee\Gc,\P;A)} \Inf_{\b \in L^0(\Omega,\Gc,\P;B)}  \E \big[ F_\varphi(t_{_0},\xi_{_0},\mu_{_0},\a,\b,\P_{(\a,\b)}) \big] \ \leq \ - 3\,\eps.
\end{equation}
We can now apply the Jankov-von Neumann measurable selection theorem, and in particular Proposition 7.50 in \cite{BertsekasShreve}, to deduce the existence of an analytically measurable function $\psi\colon L^0(\Omega,\Fc_T\vee\Gc,\P;A)\rightarrow L^0(\Omega,\Gc,\P;B)$ such that
\begin{equation}\label{psi}
\E\big[ F_\varphi(t_{_0},\xi_{_0},\mu_{_0},\a,\psi(\a),\P_{(\a,\psi(\a))}) \big] \ \leq \ - 2\,\eps, \qquad \text{for every }\a\in L^0(\Omega,\Fc,\P;A).
\end{equation}
More precisely, in order to apply Proposition 7.50 in \cite{BertsekasShreve}, we begin noting that $X$, $Y$, $D$, $f$, $f^*$ in \cite{BertsekasShreve} are given respectively by $L^0(\Omega,\Fc_T\vee\Gc,\P;A)$, $L^0(\Omega,\Gc,\P;B)$, $L^0(\Omega,\Fc_T\vee\Gc,\P;A)\times L^0(\Omega,\Gc,\P;B)$, $\E[ F_\varphi(t_{_0},\xi_{_0},\mu_{_0},\cdot,\cdot,\P_{(\cdot,\cdot)}) ]$, and
\[
f^*(\a) \ = \ \Inf_{\b \in L^0(\Omega,\Gc,\P;B)}  \E \big[ F_\varphi(t_{_0},\xi_{_0},\mu_{_0},\a,\b,\P_{(\a,\b)}) \big].
\]
We also introduce the following metric on $L^0(\Omega,\Fc_T\vee\Gc,\P;A)$:
\[
\rho_{L^0,A}(\alpha,\alpha') \ := \ \E\big[ \rho_A(\a,\a')\big],
\]
for any $\a,\a'\in L^0(\Omega,\Fc_T\vee\Gc,\P;A)$, where we recall that $\rho_A$ is a bounded metric on the Polish space $A$. Notice that the metric space $(L^0(\Omega,\Fc_T\vee\Gc,\P;A),\rho_{L^0,A})$ is complete. Moreover, proceeding as for the metric space $(\Ac^t,\rho_{\textup{Kr}})$ in the proof of Proposition \ref{P:Classical}, using that by assumption the $\sigma$-algebra $\Fc_T\vee\Gc$ is countably generated, we can prove that $(L^0(\Omega,\Fc_T\vee\Gc,\P;A),\rho_{L^0,A})$ is also separable, so $(L^0(\Omega,\Fc_T\vee\Gc,\P;A),\rho_{L^0,A})$ is a Polish space. Analogously, we introduce a metric $\rho_{L^0,B}$ on $L^0(\Omega,\Gc,\P;B)$ defined similarly to $\rho_{L^0,A}$; so, in particular, $(L^0(\Omega,\Gc,\P;B),\rho_{L^0,B})$ is also a Polish space. Then, we notice that the function $(\a,\b)\mapsto\E[ F_\varphi(t_{_0},\xi_{_0},\mu_{_0},\a,\b,\P_{(\a,\b)}) ]$ is Borel measurable, so, in particular, it is a lower semianalytic function. Then, by Proposition 7.50 in \cite{BertsekasShreve} it follows that: there exists an analytically measurable function $\psi\colon L^0(\Omega,\Fc_T\vee\Gc,\P;A)\rightarrow L^0(\Omega,\Gc,\P;B)$ such that
\[
\E\big[ F_\varphi(t_{_0},\xi_{_0},\mu_{_0},\a,\psi(\a),\P_{(\a,\psi(\a))}) \big] \ \leq \ f^*(\a) + \eps \ \leq \ - 2\,\eps, \quad \text{for every }\a\in L^0(\Omega,\Fc_T\vee\Gc,\P;A),
\]
where the second inequality follows from \eqref{psi1}. This concludes the proof of \eqref{psi}.

Now, by the continuity of the map $(t,\xi,\mu)$ $\in$ $[0,T]\times L^2(\Omega,\Fc_T\vee\Gc,\P;\R^n)\times\mathscr P_{\text{\tiny$2$}}(\R^n)$ $\mapsto$ 
$\E\big[ F_\varphi(t,\xi,\mu,\a,\b,\P_{(\a,\b)}) \big]$ uniform with respect to $\a$ $\in$ $L^0(\Omega,\Fc_T\vee\Gc,\P;A)$ and $\b$ $\in$ $L^0(\Omega,\Gc,\P;B)$, there exists 
$\delta$ $>$ $0$ such that, for all $\a \in L^0(\Omega,\Fc_T\vee\Gc,\P;A)$,
\beq \label{Fepssous}
\E \big[ F_\varphi(s,\xi,\mu,\a,\psi(\a),\P_{(\a,\psi(\a))}) \big] & \leq & - \eps, \;\;\; \forall s \in [t_{_0},t_{_0}+\delta], \; (\xi,\mu) \in B_\delta(\xi_{_0},\mu_{_0}).
\enq
As in \reff{hpetit}, we can take some $h$ $>$ $0$ small enough  so that, for all $\alpha$ $\in$ $\Ac$, $\beta$ $\in$ $\Bc$,
\beqs
\big(s,X_s^{t_{_0},\xi_{_0},\alpha,\beta},\P_{\text{\tiny$X_s^{t_{_0},\xi_{_0},\alpha,\beta}$}} \big) \in [t_{_0},t_{_0}+\delta]\times B_\delta(\xi_{_0},\mu_{_0}), \;\;\; 
\mbox{ for all } s \in [t_{_0},t_{_0} + h]. 
\enqs
Now, define
\[
\hat\beta[\alpha]_s \ := \ \psi(\alpha_s), \qquad \text{for every }\alpha\in\Ac.
\]
Notice that $\hat\beta[\cdot]\in\Bc_{\textup{\tiny str}}$. As a matter of fact, this follows from the following items:
\begin{enumerate}
\item[1)] For every $\alpha\in\Ac$, the map $\boldsymbol\alpha\colon s\mapsto\alpha_s(\cdot)$, from $[0,T]$ to $L^0(\Omega,\Fc_T\vee\Gc,\P;A)$, is Borel measurable.
\item[2)] For any analytically measurable map $\boldsymbol\beta\colon[0,T]\rightarrow L^0(\Omega,\Gc,\P;B)$, the process $\beta\colon\Omega\times[0,T]\rightarrow B$ defined as
\[
\beta_s(\cdot) \ = \ \boldsymbol\beta(s)(\cdot), \qquad \text{for every }s\in[0,T],
\]
is $\Gc\otimes\Bc([0,T])$-measurable; so, in particular, $\beta$ is $(\Fc_s\vee\Gc)_s$-progressively measurable, that is $\beta\in\Bc$.
\end{enumerate}
Suppose that 1) and 2) hold. Then, by 1) and the measurability property of $\psi$, for every $\alpha\in\Ac$ the map $s\mapsto\psi(\alpha_s)$ is analytically measurable from $[0,T]$ to $L^0(\Omega,\Gc,\P;B)$. Therefore, the fact that $\hat\beta[\cdot]\in\Bc_{\textup{\tiny str}}$ follows from item 2). Concerning the proof of item 1), notice that when $\alpha$ is a step process the result clearly holds; for a generic $\alpha$ the claim follows by an approximation argument. Similarly, regarding item 2), we begin noting that any analytically measurable map $\boldsymbol\beta\colon[0,T]\rightarrow L^0(\Omega,\Gc,\P;B)$ is in particular Lebesgue measurable, that is $\boldsymbol\beta$ is a measurable map from $([0,T],\mathscr L([0,T]))$ into $(L^0(\Omega,\Gc,\P;B),\Bc(L^0(\Omega,\Gc,\P;B)))$. If, for a moment, we replace $(L^0(\Omega,\Gc,\P;B),\Bc(L^0(\Omega,\Gc,\P;B)))$ by $(\R,\Bc(\R))$, then it is well-known that any $\boldsymbol\beta$ can be approximated by a sequence of step functions $\{\boldsymbol\beta_n\}_n$; since $L^0(\Omega,\Gc,\P;B)$ is separable, the same result holds true for a Lebesgue measurable map $\boldsymbol\beta$ taking values in $L^0(\Omega,\Gc,\P;B)$. As a consequence, it is enough to prove item 2) for $\boldsymbol\beta$ step function, since afterwards the claim follows by an approximation argument. When $\boldsymbol\beta$ is a step function, it is easy to see that the map $\boldsymbol\beta(\cdot)(\cdot)$ is $\Gc\otimes\Bc([0,T])$-measurable, which concludes the proof of item 2).
  
Now, by applying It\^o's formula to  $\varphi(r,\P_{\text{\tiny$X_{r}^{t_{_0},\xi_{_0},\alpha,\hat\beta[\alpha]}$}})$ between  $t_{_0}$ and $t_{_0}+h$, we get by \reff{Fepssous}
\begin{align*}
&\varphi\big(t_{_0}+h,\P_{\text{\tiny$X_{t_{_0}+h}^{t_{_0},\xi_{_0},\alpha,\hat\beta[\alpha]}$}}\big) -   \varphi(t_{_0},\mu_{_0}) \\
&= \ \E\bigg[\int_{t_{_0}}^{t_{_0}+h} 
\Big(F_\varphi\big(s,X_s^{t_{_0},\xi_{_0},\alpha,\hat\beta[\alpha]},\P_{\text{\tiny$X_{s}^{t_{_0},\xi_{_0},\alpha,\hat\beta[\alpha]}$}},
\alpha_s,\hat\beta[\alpha]_s,\P_{\text{\tiny$(\alpha_s,\hat\beta[\alpha]_s)$}} \big) \\
&\quad \ - f\big(X_s^{t_{_0},\xi_{_0},\alpha,\hat\beta[\alpha]},\P_{\text{\tiny$X_s^{t_{_0},\xi_{_0},\alpha,\hat\beta[\alpha]}$}},\alpha_s,\hat\beta[\alpha]_s,
\P_{\text{\tiny$(\alpha_s,\hat\beta[\alpha]_s)$}}\big)\Big)\,ds \bigg]  \\
&\leq \ - \eps\,h - \E\bigg[  \int_{t_{_0}}^{t_{_0}+h}  f\big(X_s^{t_{_0},\xi_{_0},\alpha,\hat\beta[\alpha]},\P_{\text{\tiny$X_s^{t_{_0},\xi_{_0},\alpha,\hat\beta[\alpha]}$}},\alpha_s,\hat\beta[\alpha]_s,\P_{\text{\tiny$(\alpha_s,\hat\beta[\alpha]_s)$}}\big)\,ds \bigg]. 
\end{align*}
Recalling that $\mathscr V(t_{_0},\mu_{_0})$ $=$   $\varphi(t_{_0},\mu_{_0})$, $\mathscr V$ $\leq$ $\varphi$, and that $\alpha$ is arbitrary in $\Ac$, we obtain
\begin{align*}
\mathscr V(t_{_0},\mu_{_0}) -  \eps\,h \ \geq \    
\sup_{\alpha\in\Ac}  \E\bigg[  \int_{t_{_0}}^{t_{_0}+h}  f\big(X_s^{t_{_0},\xi_{_0},\alpha,\hat\beta[\alpha]},\P_{\text{\tiny$X_s^{t_{_0},\xi_{_0},\alpha,\hat\beta[\alpha]}$}},\alpha_s,\hat\beta[\alpha]_s,\P_{\text{\tiny$(\alpha,\hat\beta[\alpha]_s)$}}\big)\,ds& \\
+ \, \mathscr V\big(t_{_0}+h,\P_{\text{\tiny$X_{t_{_0}+h}^{t_{_0},\xi_{_0},\alpha,\hat\beta[\alpha]}$}}\big) &\bigg], 
\end{align*}
which contradicts the DPP relation \reff{DPPvinv}. 
\ep

\begin{Remark} \label{remunivisco}
{\rm
(Uniqueness of viscosity solutions) 

\noindent Once we have the viscosity property of the value function to the dynamic programming Bellman-Isaacs equation, the next step is to state a  comparison principle for this PDE in order to 
get the characterization of the value function as the unique viscosity solution to the Bellman-Isaacs equation. Comparison principle for PDE in Wasserstein space of probability measures, and more generally on metric spaces, can be found e.g. in  \cite{ganswi15}, \cite{gantud17}, but  concern first-order equations, and to the best of our knowledge, the proof of a comparison principle for second-order equations as in \reff{lowerBI} (or \reff{upperBI}) remains a challenging issue.  On the other hand, there is a well developed theory of viscosity solutions for second-order equations of Bellman type related notably to stochastic control in Hilbert spaces, see the recent book \cite{fabgozswi17}.   In particular, we can use comparison principle  in Theorem 3.50 of this book, and check that the assumptions of this theorem  are satisfied in our context for the lifted Hamiltonian $\Hc_-$ (and $\Hc_+$) under {\bf (A1)}-{\bf (A2)}. Assuming that the function $h$ in {\bf (A1)}(ii)  
satisfies a polynomial growth condition, we then deduce from Theorem 3.50 in \cite{fabgozswi17} and our Theorem \ref{prop:visco} that the lower (resp. upper) value function $v$ (resp. $u$) is a
the unique viscosity solution to the lower (resp. upper) Bellman-Isaacs equation  \reff{lowerBIlift} (resp. \reff{upperBIlift}) on  $[0,T]\times L^2(\Omega,\Gc,\P;\R^n)$ satisfying a polynomial growth condition. 
}
\epR
\end{Remark}

\begin{Remark}
{\rm  If the Isaacs condition holds, that is $\Hc_-$ $=$ $\Hc_+$ (or equivalently $\mathscr H_-$ $=$ $\mathscr H_+$), then the Bellman-Isaacs equations  \reff{lowerBIlift} and \reff{upperBIlift} 
coincide, and by uniqueness of viscosity solutions to these equations (see Remark \ref{remunivisco}), the lower value function $v$  is equal to the upper value function $u$ (and then  
the inverse-lifted  lower  value function $\mathscr V$ is equal to the inverse-lifted upper value function  $\mathscr U$),  which means that the McKean-Vlasov stochastic differential game has a value. 
}
\epR
\end{Remark}


\appendix

\renewcommand\thesection{}

\section{}

\renewcommand\thesection{\Alph{subsection}}

\renewcommand\thesubsection{Appendix.}

\subsection{A technical lemma}
\label{App:Lemma}

\setcounter{Theorem}{0}
\setcounter{Definition}{0}
\setcounter{Proposition}{0}
\setcounter{Assumption}{0}
\setcounter{Lemma}{0}
\setcounter{Corollary}{0}
\setcounter{Remark}{0}
\setcounter{Example}{0}

\begin{Lemma}\label{L:AppTechn}
Consider a complete probability space $(\Omega,\Hc,\P)$, two random variables $\Gamma,\tilde{\Gamma}\colon\Omega\rightarrow[0,1]$ with uniform distribution, and a Borel measurable function $\phi\colon[0,1]\rightarrow[0,1]$. Suppose that
\[
\tilde{\Gamma} \ = \ \phi(\Gamma), \qquad \P\text{-a.s.}
\]
So, in particular, $\phi$ is a uniform distribution preserving map. Then, there exists a Borel measurable function $\rho\colon[0,1]\rightarrow[0,1]$ such that $\rho(\phi)(y)=y$, $\lambda$-a.e., hence
\[
\Gamma \ = \ \rho(\tilde{\Gamma}), \qquad \P\text{-a.s.}
\]
\end{Lemma}
\textbf{Proof.}
We split the proof into three steps, which can be summarized as follows:
\begin{itemize}
\item in \textbf{Step I} we consider the analytically (hence Lebesgue) measurable right-inverse $\psi$ of $\phi$ given by the Jankov-von Neumann measurable selection theorem; afterwards, we take a Borel-measurable version $\tilde\psi$ of $\psi$, which coincides with $\psi$ outside a Borel null-set $N$ (we need $\tilde\psi$, rather than $\psi$, in order to apply once again the Jankov-von Neumann theorem in \textbf{Step II}); we end \textbf{Step I} proving that, as expected, $\tilde\psi$ is a right-inverse of $\phi$ outside of the set $N$, that is $\phi(\tilde\psi)(y)=y$, $\lambda(dy)$-a.e.; our aim is then to prove that the claim follows with $\rho:=\tilde\psi$ (statement \eqref{rho=tildepsi}), namely that $\tilde\psi$ is also a left-inverse of $\phi$;
\item in \textbf{Step II} we apply the Jankov-von Neumann theorem in order to construct an analytically measurable right-inverse $\phi'$ of $\tilde\psi$, that is $\tilde\psi(\phi')(y)=y$, $\forall\,y\in[0,1]$; then, we notice that the claim follows if we prove that $\phi=\phi'$; finally, we show that this latter equality follows if $\phi'$ is a uniform distribution preserving map, namely \eqref{phi'} holds;
\item in \textbf{Step III} we prove that $\phi'$ is a uniform distribution preserving map using that both $\Gamma$ and $\tilde\Gamma$ has uniform distribution.
\end{itemize}

\vspace{2mm}

\noindent\textbf{Step I.} \emph{Borel measurable right-inverse of $\phi$.} Without loss of generality, we can suppose that $\phi$ is surjective (otherwise, we proceed along the same lines as for the function $\chi$ in \textbf{Step II} of the proof of Proposition \ref{P:Lift}). Now, using the Jankov-von Neumann selection Theorem (in particular, Corollary 18.23 in \cite{Selection}), we deduce the existence of a measurable function $\psi\colon([0,1],\mathscr L([0,1]))\rightarrow([0,1],\Bc([0,1]))$ satisfying:
\[
\phi(\psi(y)) \ = \ y, \; \text{for any $y\in[0,1]$}; \qquad\qquad \phi^{-1}(\psi^{-1}(\mathscr B)) \ = \ \mathscr B, \; \text{for any subset $\mathscr B$ of $[0,1]$}.
\]
It is well-known (see e.g. Exercise 14, Chapter 2, in \cite{Rudin}) that there exists a Borel measurable function $\tilde\psi\colon([0,1],\Bc([0,1]))\rightarrow([0,1],\Bc([0,1]))$ such that
\begin{equation}\label{tilde_psi}
\psi(y) \ = \ \tilde\psi(y), \qquad \lambda(dy)\text{-a.e.}
\end{equation}
So, in particular, there exists a $\lambda$-null set $N\in\Bc([0,1])$ such that $\psi(y)=\tilde\psi(y)$, for any $y\in[0,1]\backslash N$.

We conclude this step recalling the following properties of $\tilde\psi$, which can be deduced from the properties of $\psi$:
\begin{itemize}
\item[$1)$] $\phi(\tilde\psi(y))=y$, for any $y\in[0,1]\backslash N$ (that is, $\phi(\tilde\psi(y))=y$, $\lambda(dy)$-a.e.);
\item[$2)$] for any subset $\mathscr B$ of $[0,1]$ we have $\phi^{-1}(\tilde\psi^{-1}(\mathscr B))=\tilde{\mathscr B}$, for some subset $\tilde{\mathscr B}$ of $[0,1]$ such that: $\tilde{\mathscr B}\Delta\mathscr B\in\mathscr L([0,1])$ and $\tilde{\mathscr B}\Delta\mathscr B$ is a $\lambda$-null set.

\emph{Proof of item $2)$.} Fix a subset $\mathscr B$ of $[0,1]$. By \eqref{tilde_psi} we deduce that there exists a subset $N_{\mathscr B}\subset N$ (so, in particular, $N_{\mathscr B}$ belongs to $\mathscr L([0,1])$ and is a $\lambda$-null set) such that $\psi^{-1}({\mathscr B})\Delta\tilde\psi^{-1}({\mathscr B})=N_{\mathscr B}$. In other words, there exist two $\lambda$-null sets $N_{\mathscr B}',N_{\mathscr B}''\in\mathscr L([0,1])$ such that $\tilde\psi^{-1}({\mathscr B})=(\psi^{-1}({\mathscr B})\backslash N_{\mathscr B}')\cup N_{\mathscr B}''$. Hence
\begin{align*}
\phi^{-1}(\tilde\psi^{-1}({\mathscr B})) \ = \ \phi^{-1}((\psi^{-1}({\mathscr B})\backslash N_{\mathscr B}')\cup N_{\mathscr B}'') \ &= \ (\phi^{-1}(\psi^{-1}({\mathscr B}))\backslash\phi^{-1}(N_{\mathscr B}'))\cup\phi^{-1}(N_{\mathscr B}'') \\
&= \ ({\mathscr B}\backslash\phi^{-1}(N_{\mathscr B}'))\cup\phi^{-1}(N_B'') \ =: \ \tilde{\mathscr B}.
\end{align*}
It remains to prove that $\phi^{-1}(N_{\mathscr B}')$ (and similarly $\phi^{-1}(N_{\mathscr B}'')$) belongs to $\mathscr L([0,1])$ and is a $\lambda$-null set. We begin noting that since $N_{\mathscr B}'\in\mathscr L([0,1])$, by definition of $\mathscr L([0,1])$ (see e.g. Theorem 1.36 in \cite{Rudin}), there exists a $\lambda$-null set $\hat N_{\mathscr B}'\in\Bc([0,1])$ such that $N_{\mathscr B}'\subset\hat N_{\mathscr B}'$. Since $\phi^{-1}(N_{\mathscr B}')\subset\phi^{-1}(\hat N_{\mathscr B}')$, it is enough to prove that $\phi^{-1}(\hat N_{\mathscr B}')$ belongs to $\mathscr L([0,1])$ and is a $\lambda$-null set. Actually, $\phi^{-1}(\hat N_{\mathscr B}')$ belongs to $\Bc([0,1])$ (so, in particular, to $\mathscr L([0,1])$) since $\hat N_B'\in\Bc([0,1])$ and $\phi$ is a measurable function from $([0,1],\Bc([0,1]))$ into $([0,1],\Bc([0,1]))$. Now, recalling that $\lambda$ is the distribution of $\Gamma$ and $\tilde{\Gamma}$, we obtain
\[
\lambda(\phi^{-1}(\hat N_{\mathscr B}')) \ = \ \P(\Gamma\in\phi^{-1}(\hat N_{\mathscr B}')) \ = \ \P(\phi(\Gamma)\in\hat N_{\mathscr B}') \ = \ \P(\tilde{\Gamma}\in\hat N_{\mathscr B}') \ = \ \lambda(\hat N_{\mathscr B}') \ = \ 0.
\]
This concludes the proof of item $2)$. 
\end{itemize}
Finally, as for $\phi$, we can suppose $\tilde\psi$ to be surjective (this property will be used in the next step in order to apply the Jankov-von Neumann theorem).

Our aim is to prove the following:
\begin{equation}\label{rho=tildepsi}
\tilde\psi(\phi(y)) \ = \ y, \qquad \lambda(dy)\text{-a.e.}
\end{equation}
that is $\rho:=\tilde\psi$ is also a $\lambda$-a.e. left-inverse of $\phi$. 

\vspace{2mm}

\noindent\textbf{Step II.} \emph{The function $\phi'$.} We apply the Jankov-von Neumann measurable selection theorem to the function $\tilde\psi$, so we deduce the existence of a measurable function $\phi'\colon([0,1],\mathscr L([0,1]))\rightarrow([0,1],\Bc([0,1]))$ satisfying:
\[
\tilde\psi(\phi'(y)) \ = \ y, \; \text{for any $y\in[0,1]$}; \qquad\qquad \tilde\psi^{-1}((\phi')^{-1}(\mathscr B)) \ = \ \mathscr B, \; \text{for any subset $\mathscr B$ of $[0,1]$}.
\]
The claim follows (see, in particular, \eqref{rho=tildepsi}) if we prove that
\begin{equation}\label{claim}
\phi \ = \ \phi', \qquad \lambda\text{-a.e.}
\end{equation}
In order to prove \eqref{claim}, notice that
\[
\phi(\tilde\psi(\phi'(y))) \ = \ \phi(y), \qquad \text{for all }y\in[0,1]
\]
and, by property 1) above,
\[
\phi(\tilde\psi(\phi'(y))) \ = \ \phi'(y), \qquad \text{for all }\phi'(y)\in[0,1]\backslash N.
\]
Hence
\[
\phi \ = \ \phi', \qquad \text{on }\big\{y\colon\phi'(y)\notin N\big\}.
\]
It remains to prove that the set $\{y\colon\phi'(y)\notin N\}$ is a $\lambda$-null set. This holds true if we show that $\phi'$ is a uniform distribution preserving map, namely
\begin{equation}\label{phi'}
\lambda((\phi')^{-1}(\mathscr B)) \ = \ \lambda(\mathscr B), \qquad \text{for every Borel subset }\mathscr B\subset[0,1].
\end{equation}

\vspace{2mm}

\noindent\textbf{Step III.} \emph{$\phi'$ is a uniform distribution preserving map.} Suppose for a moment that both random variables
\[
\tilde\psi(\phi(\Gamma)) \qquad \text{ and } \qquad \phi'(\tilde\psi(\phi(\Gamma)))
\]
have uniform distribution. Then, for any Borel subset $\mathscr B$ of $[0,1]$,
\[
\lambda(\mathscr B) \ = \ \P(\phi'(\tilde\psi(\phi(\Gamma)))\in\mathscr B) \ = \ \P(\tilde\psi(\phi(\Gamma))\in(\phi')^{-1}(\mathscr B)) \ = \ \lambda((\phi')^{-1}(\mathscr B)),
\]
which implies that \eqref{phi'} holds true. Finally, we report below the proof that $\tilde\psi(\phi(\Gamma))$ and $\phi'(\tilde\psi(\phi(\Gamma)))$ are uniformly distributed.

\vspace{1mm}

\noindent\emph{The random variable $\tilde\psi(\phi(\Gamma))$.} Since, by assumption, $\Gamma$ has uniform distribution, $\tilde\psi(\phi(\Gamma))$ also has uniform distribution, as a matter of fact
\[
\P(\tilde\psi(\phi(\Gamma))\in\mathscr B) \ = \ \P(\Gamma\in\phi^{-1}(\tilde\psi^{-1}(\mathscr B))) \ = \ \P(\Gamma\in\tilde{\mathscr B}) \ = \ \lambda(\tilde{\mathscr B}) \ = \ \lambda(\mathscr B),
\]
where we have used the properties of $\tilde{\mathscr B}:=\phi^{-1}(\tilde\psi^{-1}(\mathscr B))$ stated in item 2) above.

\vspace{1mm}

\noindent\emph{The random variable $\phi'(\tilde\psi(\phi(\Gamma)))$.} Since, by assumption, $\tilde\Gamma=\phi(\Gamma)$ has uniform distribution, $\phi'(\tilde\psi(\phi(\Gamma)))$ also has uniform distribution, as a matter of fact
\[
\P(\phi'(\tilde\psi(\phi(\Gamma)))\in\mathscr B) \ = \ \P({\tilde{\Gamma}}\in\tilde\psi^{-1}((\phi')^{-1}(\mathscr B))) \ = \ \P({\tilde{\Gamma}}\in\mathscr B) \ = \ \lambda(\mathscr B).
\]
\ep

\vspace{9mm}

\small
\bibliographystyle{plain}
\bibliography{references}

\end{document}